\documentclass[11pt]{article}
\usepackage{amssymb}
\usepackage{amsmath}
\usepackage{amsfonts}
\newcommand{\prob}{\mathop{\mathbb{P}}}
\newcommand{\proba}[1]{\mathop{\mathbb{P}} \left( #1 \right)}
\newcommand{\esp}[1]{\mathop{\mathbb{E}} \left( #1 \right)}
\newcommand{\var}[1]{\mathop{\mathbb{\tilde{V}}} \left( #1 \right)}
\newcommand{\R}{\mathop{\mathbb{R}}}
\newcommand{\N}{\mathop{\mathbb{N}}}
\newcommand{\ind}{\mathop{\mathrm{1 \! I}}}

\newcommand{\AF}{strong approximation }
\newcommand{\AFS}{strong approximations }
\newcommand{\F}{{\cal F}}

\newcommand{\T}{\tilde{C}_1}
\newcommand{\TT}{\tilde{\tilde{C}}_1}

\newcommand{\DA}{ {\Delta_A}^{i,l(i,v)}_{j,k(j,u)} }
\newcommand{\DB}{ {\Delta_B}^{i,l(i,v)}_{j,k(j,u)} }
\newcommand{\DC}{ {\Delta_C}^{i,l(i,v)}_{j,k(j,u)} }
\newcommand{\DDC}{ {C}^{i,l(i,v)} }
\newcommand{\X}{{\xi}^{i,l(i,v)}_{j,k(j,u)}}
\newcommand{\SI}{\left( \sigma^{i,l(i,v)}_{j,k(j,u)} \right)^2}

\newtheorem{defi}{Definition}[section]
\newtheorem{theo}{Theorem}[section]
\newtheorem{lem}[theo]{Lemma}
\newtheorem{ineg}[theo]{Inequalities}
\numberwithin{equation}{section}
\title{Improvement of two Hungarian bivariate theorems.}
\author{CASTELLE Nathalie \\
Laboratoire de Math\'ematiques - UMR 8628 \\ 
B\^at. 425, Universit\'e
de Paris-Sud \\91405 Orsay Cedex, FRANCE \\
e-mail: Nathalie.Castelle@math.u-psud.fr}
\begin{document}
\maketitle
\noindent
\begin{center} {\bf R\'esum\'e} \end{center}
Nous introduisons une nouvelle technique pour \'etablir des
th\'eor\`emes hongrois multivari\'es. Appliqu\'ee dans 
cet article aux th\'eor\`emes bivari\'es d'approximation forte du processus 
empirique uniforme, cette technique am\'eliore le
r\'esultat de Koml\'os, Major et Tusn\'ady (1975) ainsi que
les n\^otres (1998). Plus pr\'ecis\'ement,
nous montrons que l'erreur dans l'approximation du $n$-processus
empirique uniforme bivari\'e par un pont brownien bivari\'e
est d'ordre $n^{-1/2}(\log (nab))^{3/2}$ sur le pav\'e 
$[0,a]\times[0,b]$, $0\leq a, b \leq 1$, 
et que l'erreur dans l'approximation du
$n$-processus empirique uniforme univari\'e par un processus de Kiefer
est d'ordre $n^{-1/2}(\log (na))^{3/2}$ sur 
l'intervalle $[0,a]$, $0\leq a \leq 1$.
Dans les deux cas la borne d'erreur globale est donc d'ordre 
$n^{-1/2}(\log (n))^{3/2}$. Les r\'esultats pr\'ec\'edents
donnaient depuis l'article de 1975 de Koml\'os, Major et Tusn\'ady 
une borne d'erreur globale d'ordre $n^{-1/2}(\log (n))^{2}$,
et depuis notre article de 1998 des bornes d'erreur locales d'ordre
$n^{-1/2}(\log (nab))^{2}$ ou $n^{-1/2}(\log (na))^{2}$. 
Le nouvel argument de cet article consiste \`a reconna\^\i tre des 
martingales dans les termes d'erreur, puis \`a leur appliquer une 
in\'egalit\'e exponentielle de Van de Geer (1995) ou 
de de la Pe\~na (1999).
L'id\'ee est de borner le compensateur du terme d'erreur,
au lieu de borner le terme d'erreur lui-m\^eme.
\begin{center} {\bf Abstract} \end{center}
We introduce a new technique to establish Hungarian multivariate theorems.
In this article, we apply this technique
to the strong approximation bivariate theorems 
of the uniform empirical process. It improves the Komlos, Major and 
Tusn\'ady's result (1975) as well as our own (1998). More precisely, we show 
that the error in the approximation of the uniform bivariate
$n$-empirical process by a bivariate Brownian bridge is of order
$n^{-1/2}(\log (nab) )^{3/2}$ on the rectangle $[0,a]\times[0,b]$,
$0\leq a,b \leq 1$,  
and that the error in the approximation of the uniform univariate
$n$-empirical process by a Kiefer process is of order 
$n^{-1/2}(\log (na) )^{3/2}$ on the interval $[0,a]$, $0\leq a \leq 1$. 
In both cases, the global error bound is therefore of order
$n^{-1/2}(\log (n) )^{3/2}$. Previously, from the 1975 article of
Komlos, Major and Tusn\'ady, 
the global error bound was of order 
$n^{-1/2}(\log (n) )^{2}$,
and from our 1998 article, the local error bounds were
of order $n^{-1/2}(\log (nab) )^{2}$ or 
$n^{-1/2}(\log (na) )^{2}$. The new feature of this 
article is to identify martingales in the error terms and to apply to them
a Van de Geer's (1995) or de la Pe\~na's (1999) exponential inequality. The idea
is to bound of the compensator of the error term,
instead of bounding of the error term itself.
\\ \\ {\bf AMS 1991: 60F17, 60G15, 60G42, 62G30}. 
\section{Introduction and results.}
Let $(X_i,Y_i),i\geq 1$ be a sequence of independent and identically 
distribued random couples with uniform on $[0,1]^2$ distribution, defined on
the same probability space $(\Omega,{\cal A},\prob)$.
We assume that $\Omega$ is rich enough so that there exists on 
$(\Omega,{\cal A},\prob)$ a variable with uniform distribution on
$[0,1]$ independent of the sequence $(X_i,Y_i),i\geq 1$.
Let us denote by \( {H_n} \) the cumulative empirical distribution function
associated with $(X_i,Y_i),i=1,\ldots,n$:
\[ {H}_n(t,s)=\frac{1}{n} \sum_{i=1}^n \ind{_{X_i\leq t,Y_i\leq s}} \]
   for $ (t,s)\in [0,1]^2 $,
and let us denote by ${F_n}$, ${G_n}$ the univariate cumulative empirical  
distribution functions : ${F_n}(t)={H_n}(t,1), \; {G_n}(s)={H_n}(1,s)$.
Let us recall the definitions of the Gaussian processes which appear in the strong
approximation theorems of these cumulative empirical distribution functions.
\begin{defi} A Brownian bridge $B$ is a continuous Gaussian process 
defined on $[0,1]$ such that 
$\esp{B(t)}=0$, $\esp{B(t)B(t')}=(t\wedge t')-tt'$.
A bivariate Brownian bridge $D$ is a continuous Gaussian process defined on
$[0,1]^2$ such that
$\esp{D(t,s)}=0$, $\esp{D(t,s)D(t',s')}=(t\wedge t')(s\wedge s')-tt'ss'$.

\end{defi}
\begin{defi} A Kiefer process $K$ is a continuous Gaussian
process defined on $[0,1]\times [0,1]$ such that
$\esp{K(t,s)}=0$, $\esp{K(t,s)K(t',s')}=(s\wedge s')((t\wedge t')-tt')$.
We call Kiefer process on $[0,1]\times \N$ or on
$[0,1]\times \{0,\ldots,n\}$ a Gaussian process such that
$\esp{K(t,j)}=0$, $\esp{K(t,j)K(t',j')}=(j\wedge j')((t\wedge t')-tt')$.
In this case, $K$ may be defined as a sum of independent Brownian bridges :
$K(t,0)=0$, \( K(t,j)=\sum_{i=1}^j B_i(t) \).
\end{defi}
In their famous paper of 1975,
Koml\'os, Major et Tusn\'ady established the \AF 
of the univariate cumulative empirical distribution function 
by a Brownian bridge, and also by a Kiefer process.
This last approximation, more powerful, was in fact a bivariate
approximation. The paper of 1975 left many questions open. After wards, were
carried out the \AF of the bivariate cumulative empirical distribution function 
(Tusn\'ady (1977a), Castelle et Laurent (1998)), and also
the univariate local \AF (Mason et Van Zwet (1987))
and the bivariate local \AFS (Castelle et Laurent (1998)). 
These results are summarized by the two following theorems
(Castelle (2002)). In these theorems, and throughout this article, 
we denote by log the function $x\rightarrow \ln (x \vee e)$.
\begin{theo} \label{th1} 
Let ${H}_n$ be the bivariate cumulative empirical distribution function
previously defined. For any integer $n$, there exists a bivariate Brownian
bridge $D^{(n)}$ such that for all positive $x$ and for all
$(a,b)\in [0,1]^2$ we have :
\begin{eqnarray} \label{in1} \proba { \sup_{0\leq t \leq a,0\leq s\leq b}
   |n{H_n}(t,s) -nts - \sqrt{n} {D}^{(n)}(t,s)|
    \geq  (x+C_1 \log (nab)) \log (nab) } \leq \Lambda_1 
    \exp(-\lambda_1 x)  \\ \label{in1bis}
    \proba { \sup_{0\leq t \leq a}
    |n F_n(t) -nt - \sqrt{n} {D}^{(n)}(t,1)|
    \geq x+C_0 \log (na) } \leq \Lambda_0
    \exp(-\lambda_0 x)  \\ \label{in1ter}
    \proba { \sup_{0\leq s \leq b}
    |n G_n(s) -ns - \sqrt{n} {D}^{(n)}(1,s)|
    \geq x+C_0 \log (nb) } \leq \Lambda_0
    \exp(-\lambda_0 x)
\end{eqnarray}
where \( C_0,\Lambda_0,\lambda_0,C_1,\Lambda_1,\lambda_1 \) are absolute positive
constants.
\end{theo}
Remark : in the cases $a=1$ and  $b=1$, Bretagnolle and Massart (1989) proved 
Inegalities (\ref{in1bis}), (\ref{in1ter})
with $C_0=12$, $\Lambda_0=2$ and $\lambda_0=1/6$. 
\\
\begin{theo} \label{th2} 
Let $(F_j),j\geq 1$ be the sequence of univariate cumulative empirical 
distribution functions previously defined.
There exists a Kiefer process $K$ defined on $[0,1]\times \N$
such that for all positive $x$ and for all
$a\in [0,1]$ we have :
\begin{equation*}
   \proba { \sup_{1\leq j \leq n} \sup_{0\leq t \leq a}
      |j{F_j}(t) -jt -{K}(t,j)| \geq 
      (x+C_2 \log (na)) \log (na) }
      \leq \Lambda_2 \exp(-\lambda_2 x)  
\end{equation*}     
where \( C_2,\Lambda_2,\lambda_2 \) are absolute positive
constants.
\end{theo}
The questions which remain are the optimality of the error bound
in dimension 2 and the one, more general, of the \AFS of the uniform on
$[0,1]^d$, $d\geq 3$, empirical process. We think, but it is still to be proved,
in dimension $d$ the error bound for the global \AF is of order
$(\log (n))^{{(d+1)}/{2}}$, and the error bound for the local \AF 
on $[0,a_1]\times \cdots\times [0,a_d]$ is of order
$(\log (na_1\cdots a_d))^{{(d+1)}/{2}}$.
In this paper, we improve the error bound in dimension 2 and we obtain the
following results :
\begin{theo} \label{th3}
In Theorem \ref{th1} we have also the inequality
\begin{equation} \label{in3}
    \proba { \sup_{0\leq t \leq a,0\leq s\leq b}
    |n{H_n}(t,s) -nts - \sqrt{n} {D}^{(n)}(t,s)|
    \geq  (x+C_1 \log (nab))^{3/2} } \leq \Lambda_1 
    \exp(-\lambda_1 x). 
\end{equation}
\end{theo}  
\begin{theo} \label{th4}
In Theorem \ref{th2} we have also the inequality
\[
\proba { \sup_{1\leq j \leq n} \sup_{0\leq t \leq a}
      |j{F_j}(t) -jt -{K}(t,j)| \geq 
      (x+C_2 \log (na))^{3/2} }
      \leq \Lambda_2 \exp(-\lambda_2 x).  
\]
\end{theo}  
We refer now to the paper of Castelle (2002)
which establishes that Theorem \ref{th3} leads to
Theorem \ref{th4}. More precisely, Theorem \ref{th3}
is equivalent to the following theorem :
\begin{theo} \label{th5} 
Let $(F_j),j\geq 1$ be the sequence of univariate cumulative empirical 
distribution functions previously defined. For any integer $n$,
there exists a Kiefer process $K^{(n)}$  defined on $[0,1]\times \{1,\ldots,n\}$ 
such that for all positive $x$, for all
$a\in [0,1]$ and for all integer $m\leq n$ we have :
\begin{eqnarray*} 
\proba { \sup_{1\leq j \leq m} \sup_{0\leq t \leq a}
      |j{F_j}(t) -jt -{K}^{(n)}(t,j)| \geq 
      (x+C \log (ma)) \log(ma)  }
      \leq \Lambda \exp(-\lambda x)  
\\
   \proba { \sup_{1\leq j \leq m} \sup_{0\leq t \leq a}
      |j{F_j}(t) -jt -{K}^{(n)}(t,j)| \geq 
      (x+C \log (ma))^{3/2}  }
      \leq \Lambda \exp(-\lambda x)  
\\ 
\proba { \sup_{0\leq t \leq a}
      |n{F_n}(t) -nt -{K}^{(n)}(t,n)| \geq 
      x+C_0 \log (na)  }
      \leq \Lambda_0 \exp(-\lambda_0 x)  
\end{eqnarray*}      
where \( C_0,\Lambda_0,\lambda_0\) are the constants of Theorem \ref{th1}
and where \( C,\Lambda,\lambda \) are absolute positive constants.
\end{theo}
This last theorem leads easily to Theorem
\ref{th4}. Thus, the purpose of all the subsequent sections of this paper
will be dedicated to prove Theorem \ref{th3}.
\section{Construction.} \label{2}
We use the Koml\'os, Major et Tusn\'ady construction (1975). 
More expansive explanations could be found in their article, and
also in Castelle and Laurent article (1998).
It is easier to construct the empirical process on the Gaussian process than
to construct the Gaussian process on the empirical process.
Therefore we posit a bivariate Brownian bridge $D$ and we construct
$H_n$ such that Inequalities (\ref{in1}), 
(\ref{in1bis}), (\ref{in1ter}), (\ref{in3}) hold.
In this way we obtain the reversed form of Theorem \ref{th3}.
Invoking Skorohod (1976) the theorem itself works. 
\subsection{Definition of Gaussian variables used in the construction.}
If the probability space is rich enough (if there exists on 
$(\Omega,{\cal A},\prob)$ a variable with uniform distribution on
$[0,1]$ independent of $D$), there then exists a bivariate Wiener process $W$
such that \[ D(t,s)=W(t,s)-tsW(1,1). \] 
Let us denote by $W(]t_1,t_2],]s_1,s_2])$ the expression
\[ W(t_2,s_2)-W(t_1,s_2)-W(t_2,s_1)+W(t_1,s_1). \]
Let $N$ be the integer such that $2^{N-1}< n \leq 2^N$. We set
\[ Z^{i,l}_{j,k}=\sqrt{n} W\left(]\frac{k2^j}{2^N},\frac{(k+1)2^j}{2^N}],
                        ]\frac{l2^i}{2^N},\frac{(l+1)2^i}{2^N}] \right)
\]
with \( i \in \{0,\ldots,N\}, \; l \in \{0,\ldots,2^{N-i}-1 \}, \;
j\in \{0,\dots,N \}, \; k \in \{0, \ldots, 2^{N-j}-1 \} \). 
We define now a filtration.
\[
\F^{N}_{j}=\sigma \left ( Z^{N,0}_{j,k}; k \in \{0, \ldots, 2^{N-j}-1 \} \right) 
\]
and for $i<N$,
\[
\F^{i}_{j}=\sigma \left( 
\begin{array}{l} Z^{i+1,l}_{0,k}; l \in \{0,\ldots,2^{N-(i+1)}-1 \} ;
k \in \{0, \ldots, 2^{N}-1 \} \\ Z^{i,l}_{j,k}; l \in \{0,\ldots,2^{N-i}-1 \}
; k \in \{0, \ldots, 2^{N-j}-1 \} \end{array} \right).
\] 
We have
\[ \F^{i_1}_{j_1} \subset \F^{i_2}_{j_2} \mbox{ if and only if }
\left\{ \begin{array}{l} i_1>i_2 \\ \mbox{ or } \\ i_1=i_2 \mbox{ and }
j_1>j_2. \end{array} \right.  \]
In other words,
\[ \F^N_N \subset \F^N_{N-1}\subset \cdots \F^N_{0} \subset \F^{N-1}_N \subset
\cdots \subset \F^0_0 .\]
The variables used in the construction are the variables
\[ V^{i,2l}_{j,2k}=
                  \left\{ \begin{array}{ll}
                   Z^{i,2l}_{j,2k}-\esp { Z^{i,2l}_{j,2k}/\F^i_{j+1} }, 
                   & \mbox{ if } i\leq N \mbox{ and } j<N, \\
                   Z^{i,2l}_{j,2k}-\esp { Z^{i,2l}_{j,2k}/\F^{i+1}_{0} }
                   & \mbox{ if } i< N \mbox{ and } j=N.
                   \end{array} \right.
\]                   
One easily obtains 
\begin{eqnarray*}
V^{N,0}_{j,2k}=\frac{Z^{N,0}_{j,2k}-Z^{N,0}_{j,2k+1}}{2}, \\ 
V_{N,0}^{i,2l}=\frac{Z_{N,0}^{i,2L}-Z_{N,0}^{i,2l+1}}{2}, \\
V^{i,2l}_{j,2k}=\frac{1}{4}(Z^{i,2l}_{j,2k}-Z^{i,2l+1}_{j,2k}-Z^{i,2l}_{j,2k+1}
+Z^{i,2l+1}_{j,2k+1}) \mbox{ if } i<N \mbox{ and } j<N.
\end{eqnarray*} 
These variables are independent Gaussian random variables, with expectation 0
and with variance 
\begin{eqnarray*}
\mbox{Var}(V^{N,0}_{j,2k})=\frac{\gamma 2^{j}}{2}, \\
\mbox{Var}(V_{N,0}^{i,2l})=\frac{\gamma 2^{i}}{2}, \\
\mbox{Var}(V^{i,2l}_{j,2k})= \frac{\gamma 2^{i+j-N}}{4} \mbox{ if } i<N 
\mbox{ and } j<N,
\end{eqnarray*}
with $\gamma=n/2^N$.
\subsection{Construction of the empirical process.} 
Define the inverse of a function $f$ by 
\( f^{-1}(v)= \inf \{ u/ f(u) \geq v \}
\). Denote by $\Phi,\Psi_n,\Phi_{n,n_1,n_2}$ the cumulative repartition functions
of the standard normal distribution, of the binomial distribution ${\cal B}(n,1/2)$, 
of the hypergeometric distribution ${\cal H}(n,n_1,n_2)$:
\begin{eqnarray*}
\Phi(u)=\frac{1}{\sqrt{2\pi}}\int_{-\infty}^u \exp (-t^2/2) dt, \\
\Psi_n(u)=\sum_{k=0}^{[u]} \left( \begin{array}{c} n \\ k \end{array} \right)
(\frac{1}{2})^n \mbox{ for } u\in [0,n], \\
\Phi_{n,n_1,n_2}(u)=\sum_{k=0}^{[u]}
\frac{ \left( \begin{array}{c} n_2 \\ k \end{array} \right)
\left( \begin{array}{c} n-n_2 \\ n_1-k \end{array} \right) }
{ \left( \begin{array}{c} n \\ n_1 \end{array} \right) }
\mbox{ for } u \in [\max(0,n_1+n_2-n),\min(n_1,n_2)].
\end{eqnarray*}
We construct the new variables as follows :
\[ ({\cal C}_1) \left\{ \begin{array}{l} 
U^{N,0}_{N,0}=n  \\ \\
U^{N,0}_{j,2k}=\Psi^{-1}_{U^{N,0}_{j+1,k}} \circ \Phi 
( \left( \frac{\gamma 2^{j}}{2} \right)^{-1/2} V^{N,0}_{j,2k}) \\ \\
U^{N,0}_{j,2k+1}
= U^{N,0}_{j+1,k}-U^{N,0}_{j,2k}
\end{array}  \right. \]
for \(j=N-1,\ldots,0 \) and \( k\in \{ 0,\ldots,2^{N-(j+1)}-1 \} \),
\[ ({\cal C}_2) \left\{ \begin{array}{l} 
U_{N,0}^{i,2l}=\Psi^{-1}_{U_{N,0}^{i+1,l}} \circ \Phi 
( \left( \frac{\gamma 2^{i}}{2} \right)^{-1/2} V_{N,0}^{i,2l}) \\ \\
U_{N,0}^{i,2l+1}
= U^{N,0}_{j+1,k}-U^{N,0}_{j,2k}
\end{array}  \right. \]
for \(i=N-1,\ldots,0 \) and \( l\in \{ 0,\ldots,2^{N-(i+1)}-1 \} \),
\[ ({\cal C}_3) \left\{ \begin{array}{l}
U^{i,2l}_{j,2k}
=\Phi^{-1}_{U^{i+1,l}_{j+1,k},
U^{i+1,l}_{j,2k},
U^{i,2l}_{j+1,k}} \circ \Phi 
( \left( \frac{\gamma 2^{i+j-N}}{4} \right)^{-1/2} V^{i,2l}_{j,2k}) \\ \\
U^{i,2l}_{j,2k+1}=
U^{i,2l}_{j+1,k}-U^{i,2l}_{j,2k}
\\ \\ 
U^{i,2l+1}_{j,2k}=
U^{i+1,l}_{j,2k}-U^{i,2l}_{j,2k}
\\ \\ 
U^{i,2l+1}_{j,2k+1}=
U^{i+1,l}_{j+1,k}-U^{i+1,l}_{j,2k}-U^{i,2l}_{j+1,k}+U^{i,2l}_{j,2k} 
\end{array}     \right.  \]
for  \(i=N-1,\ldots,0 \); \( l\in \{ 0,\ldots,2^{N-(i+1)}-1 \} \);
\(j=N-1,\ldots,0 \) and \( k\in \{ 0,\ldots,2^{N-(j+1)}-1 \} \).
In this way, we obtain a $\R^{2^N} \otimes \R^{2^N}$ vector,
denoted by $M$, defined by
\begin{eqnarray} \nonumber
M=(U^{0,0}_{0,0},U^{0,0}_{0,1},\ldots,U^{0,0}_{0,2^N-1},
   U^{0,1}_{0,0},U^{0,1}_{0,1},\ldots,U^{0,1}_{0,2^N-1}, \cdots, \\
   U^{0,2^N-1}_{0,0},U^{0,2^N-1}_{0,1},\ldots,U^{0,2^N-1}_{0,2^N-1}). 
\label{disemp}
\end{eqnarray}
>From Proposition 3.2 of Castelle and Laurent (1998), the vector $M$ has the
multinomial distribution
\[ {\cal M}_{2^N\times 2^N}(n,(\frac{1}{2^N})^2,\ldots,(\frac{1}{2^N})^2). \] 
Remark: Proposition 3.2 of Castelle and Laurent (1998) contains two
Equalities called (3.6) and (3.7). 
The restriction {\em n even} at the beginning
of the proposition concerns only equality (3.7). In this paper, 
we use only Equality (3.6) which is valid for all integer $n$. \\ \\
Thus the vector $M$ has the same distribution as the discretization
of a $n$-empirical cumulative distribution function
on small slabs with size $\frac{1}{2^N}\times \frac{1}{2^N}$.
If there exists on $(\Omega,{\cal A},\prob)$ a variable with uniform distribution 
on $[0,1]$ independent of $W$,
Skohorod's Theorem (1976) ensures the existence
of a bivariate $n$-empirical cumulative distribution function,
which we denote by $H_n$ from now on, such that :
\[ nH_n \left(]\frac{k}{2^N},\frac{(k+1)}{2^N}],
                        ]\frac{l}{2^N},\frac{(l+1)}{2^N}] \right)
                        =U^{0,l}_{0,k} 
\]
for $l\in \{ 0,\ldots,2^{N}-1 \}$,$k\in \{ 0,\ldots,2^{N}-1 \}$.  
\subsection{Hypergeometric Lemma.}
The control of the distance between the empirical and the Gaussian processes
needs the control of the difference between the variables
$U^{i,2l}_{j,2k}$ and $V^{i,2l}_{j,2k}$. 
For steps $({\cal C}_1)$, $({\cal C}_2)$, this control is given 
by Tusn\'ady's Lemma (1977b) proved in 1989
by Bretagnolle and Massart. We don't use this part of Tusn\'ady's Lemma 
in this paper, instead  we use 
Inequalities (\ref{in1bis}), (\ref{in1ter}) 
which were proved from this lemma.
For step $({\cal C}_3)$, the control is given by a lemma, the so-called
hypergeometric Lemma, proved in 1998 by Castelle and Laurent.
\begin{lem} \label{hyper}
For all indexes $i,j \leq N-1$, we set 
\[  \delta^{i+1,l}_{j,2k}=
                   \frac{U^{i+1,l}_{j,2k}-U^{i+1,l}_{j,2k+1}}{U^{i+1,l}_{j+1,k}}                   
                   \;\;\; \mbox{ and } \;\;\; 
                   \tilde{\delta}^{i,2l}_{j+1,k}=
                   \frac{U^{i,2l}_{j+1,k}-U^{i,2l+1}_{j+1,k}}{U^{i+1,l}_{j+1,k}}.
\]
If \( |\delta^{i,2l}_{j+1,k} \tilde{\delta}^{i,2l}_{j+1,k}| \leq \epsilon^2 <1 \)
we have
\begin{eqnarray*}
|U^{i,2l}_{j,2k}-\esp { U^{i,2l}_{j,2k}/{\cal F}^i_{j+1} }
-\left( \var { U^{i,2l}_{j,2k}/{\cal F}^i_{j+1} } \right)^{1/2}
\left( \frac{ \gamma 2^{i+j-N} }{4} \right)^{-1/2} V^{i,2l}_{j,2k}| \\
\leq \alpha + \beta \left(
\left( \frac{\gamma 2^{i+j-N}}{4} \right)^{-1/2}  V^{i,2l}_{j,2k}
\right)^2. 
\end{eqnarray*}
with
\begin{eqnarray*}
\esp { U^{i,2l}_{j,2k}/{\cal F}^i_{j+1} }={U^{i+1,l}_{j+1,k}}
\frac{U^{i+1,l}_{j,2k}}{U^{i+1,l}_{j+1,k}} 
\frac{U^{i,2l}_{j+1,k}}{U^{i+1,l}_{j+1,k}}, \\
\var { U^{i,2l}_{j,2k}/{\cal F}^i_{j+1} }=
{U^{i+1,l}_{j+1,k}}
\frac{U^{i+1,l}_{j,2k}}{U^{i+1,l}_{j+1,k}} 
\frac{U^{i+1,l}_{j,2k+1}}{U^{i+1,l}_{j+1,k}} 
\frac{U^{i,2l}_{j+1,k}}{U^{i+1,l}_{j+1,k}}
\frac{U^{i,2l+1}_{j+1,k}}{U^{i+1,l}_{j+1,k}}.
\end{eqnarray*}
where $\alpha$ and $\beta$ are positive constants which depend only on
$\epsilon$. Moreover if $\epsilon^2 =1/8$ and if the condition
\( \var { U^{i,2l}_{j,2k}/{\cal F}^i_{j+1} } \geq 4.5$ holds, the constants
$\alpha =3$ and $\beta =0.41$ are appropriate.
\end{lem}
\section{Control of the approximation error.} \label{3}
Inequalities (\ref{in1}), (\ref{in1bis}) and (\ref{in1ter}) have already
been proved. Let 
$P$ be the probability to be controlled to obtain (\ref{in3}) :
\[ P=\proba { \sup_{0\leq t \leq a,0\leq s\leq b}
    |n{H_n}(t,s) -nts - \sqrt{n} {D} (t,s)|
    \geq  (x+C_1 \log (nab))^{3/2} } .  \]
Let $\T$ be a positive constant, not fixed for the moment,
but such that $\T\geq 10$. We do not try to optimize the constants
in this paper as a numeric work will be realised later.
Set
\[ C_1=\left\{ \begin{array}{ll} \displaystyle{ \frac{3 \T }{2} }+ 
\displaystyle{
       2 \left( \frac{3 C_0 }{4} \right)^{2/3} }& \mbox{ when } ab=1 \\
\displaystyle{
       \frac{3 \T }{2} + 
        2 \left( 3 C_0 \right)^{2/3} }& \mbox{ when } ab<1       
\end{array} \right. \]
where $C_0$ is the constant of Inequalities (\ref{in1bis}), 
(\ref{in1ter}).
In the case $(x/2)+\T \log (nab) \geq \gamma^2 (nab)/8$, 
the result stems not from the construction, but from maximal
inequalities for the bivariate Brownian bridge and the 
bivariate $n$-empirical bridge.
These inequalities, summarized in Inequalities \ref{maxi} below,
are due to Adler and Brown (1986), 
Talagrand (1994) and Castelle and Laurent (1998). 
\begin{ineg} \label{maxi}
a) For all $a,b\in[0,1]$ such that \( 0\leq ab \leq
1/2\) we have:
\[ \proba { \sup_{(s,t)\in[0,b]\times[0,a]} |
n({H}_n(s,t)-st)|\geq x } \leq 2e \exp(-nab(1-ab)h(\frac{x}{nab}))  \]
where the function $h$ is defined for $t>-1$ by \( h(t)=(1+t)\ln (1+t)-t \).
\\ \\
b) There exists an universal positive constant $C$
such that:  
\[ \proba { \sup_{(s,t)\in[0,1]\times[0,1]} |
   \sqrt{n}({H}_n(s,t)-st)| \geq x } \leq Cx^2 \exp(-2x^2) . \] 
c) For all $a,b\in[0,1]$ such that \( 0\leq ab \leq
1/2\) we have:
\[ \proba { \sup_{(s,t)\in[0,b]\times[0,a]} |D(s,t)|\geq x }
\leq 2e \exp(-\frac{x^2(1-ab)}{2ab}) . \]
d) There exists an universal positive constant $C$ such that:
\[ \proba { \sup_{(s,t)\in[0,1]\times[0,1]} |D(s,t)|\geq x }
\leq Cx^2\exp (-2x^2) . \]
\end{ineg}
We now consider the case $(x/2)+\T \log (nab) < \gamma^2 (nab)/8$.
In this case, we have $nab > 496$. Let $A$ and $B$ be the integers defined by
\[ {2^{A-1}} < na \leq {2^{A}} \mbox{ and }
   {2^{B-1}}< nb \leq {2^{B}}. \]
We have $8 \leq A,B \leq N$ and $2^{A+B-N} <4 (nab)$.
We discretize the variable $t$ on a grid with size
$\displaystyle{\frac{2^{A^*}}{2^N}}$ where $A^*$ is the integer defined by
\[ {2^{A^*+B-N}} < \frac{ 4 ((x/2)+\T \log (nab))}{\gamma} \leq {2^{A^*+B-N+1}} , \]
then we discretize the variable $s$ on a grid with size
$\displaystyle{\frac{2^{B^*}}{2^N}}$ where $B^*$ is the integer defined by
\[ {2^{A+B^*-N}} < \frac{ 4 ((x/2)+\T \log (nab))}{\gamma} \leq {2^{A+B^*-N+1}}. \]
We have $A+B^*=A^*+B$, $A^* \leq A-2$, $2^{A-A^*}<(nab)/31$.
Let us denote $\Delta^E_n(t,s)$ for \( nH_n(t,s)-nts \) and
let us denote $\Delta^G_n(t,s)$ for $\sqrt{n} D(t,s)$. 
Using the stationarity properties of the increments
\[ \{ \Delta^F_n(t,s)- \Delta^F_n(\alpha,s); 
\alpha\leq t \leq \beta; 0\leq s \leq s_0 \} \stackrel{{\cal D}}{=}
\{ \Delta^F_n(t,s) ; 0\leq t \leq \beta-\alpha; 0\leq s \leq s_0 \} \]
and
\[ \{ \Delta^F_n(t,s)- \Delta^F_n(t,\alpha); 
0 \leq t \leq t_0; \alpha \leq s \leq \beta \} \stackrel{{\cal D}}{=}
\{ \Delta^F_n(t,s) ; 0\leq t \leq t_0 ; 0\leq s \leq \beta-\alpha \} \]
where $F\in \{ E,G \}$, one gets, setting $m=(x+C_1 \log(nab))^{3/2}$,
\begin{eqnarray*}
P\leq {2^{A-A^*}} \proba{ \sup_{t\in [0,\frac{2^{A^*}}{2^N}]}
\sup_{s \in [0,b]} | \Delta^G_n(t,s) | 
\geq 0.1 m  } \\
+{2^{A-A^*}} \proba{ \sup_{t\in [0,\frac{2^{A^*}}{2^N}]}
\sup_{s \in [0,b]} | \Delta^E_n(t,s) | 
\geq 0.1 m  } \\
+{2^{B-B^*}} \proba{ \sup_{t\in [0,\frac{2^{A}}{2^N}]}
\sup_{s \in [0,\frac{2^{B^*}}{2^N}]} | \Delta^G_n(t,s) | 
\geq 0.1 m  } \\
+{2^{B-B^*}} \proba{ \sup_{t\in [0,\frac{2^{A}}{2^N}]}
\sup_{s \in [0,\frac{2^{B^*}}{2^N}]} | \Delta^E_n(t,s) | 
\geq 0.1 m  } \\
+ \proba{ \max_{1\leq u \leq 2^{A-A^*}} \max_{1\leq v \leq 2^{B-B^*}} 
|nH_n(\frac{u2^{A^*}}{2^N},\frac{v2^{B^*}}{2^N})-
n\frac{u2^{A^*}}{2^N}\frac{v2^{B^*}}{2^N}
-\sqrt{n} D(\frac{u2^{A^*}}{2^N},\frac{v2^{B^*}}{2^N})| \geq  
0.6 m }.
\end{eqnarray*}
The four first terms are controled by Inequalities
\ref{maxi} a) and c). To achieve the proof of Theorem
\ref{th3}, the following lemma remains to be proved :
\begin{lem} \label{principal} 
In the case $nab > 496$, we have
\begin{eqnarray*} \mathop{\mathbb{P}}
\left( \max_{1\leq u \leq 2^{A-A^*}} \max_{1\leq v \leq 2^{B-B^*}} 
|nH_n( \frac{u2^{A^*}}{2^N},\frac{v2^{B^*}}{2^N})-
n\frac{u2^{A^*}}{2^N} \frac{v2^{B^*}}{2^N}
-\sqrt{n} D (\frac{u2^{A^*}}{2^N}, \frac{v2^{B^*}}{2^N})| \right. 
\\ \geq  
\left. 0.6(x+C_1 \log(nab))^{3/2} \right)  \leq \Lambda_3 \exp (-\lambda_3 x) 
\end{eqnarray*}
where $\Lambda_3$, $\lambda_3$ are absolute positive constants.
\end{lem}
\section{Proof of Lemma \ref{principal}.} \label{4}
A subset of indexes 
\( \{ i_1,\ldots,i_d \} \) of \( \{ 1,\ldots,2^N \} \) can be 
identified with the $\R^{2^N}$ vector \( (x_1,\ldots, x_d) \) defined by :
\[ \begin{array}{l} x_i=1 \mbox{ for } i \in \{ i_1,\ldots,i_d \} \\
                    x_i=0 \mbox{ for } i \notin \{ i_1,\ldots,i_d \}.
\end{array} \]
Let us denote by $\gamma_u$ and $\delta_v$ the $\R^{2^N}$ vector
associated with \( \{ 1,\ldots, u 2^{A^*} \} \) and \( \{ 1,\ldots, v 2^{B^*} \} \).
Let $e^N_0$ be the $\R^{2^N}$ vector
associated with \( \{ 1,\ldots,  2^{N} \} \).
With these notations we have
\begin{eqnarray*} nH_n(\frac{u2^{A^*}}{2^N},\frac{v2^{B^*}}{2^N})-
   n\frac{u2^{A^*}}{2^N}\frac{v2^{B^*}}{2^N}
   -\sqrt{n} D(\frac{u2^{A^*}}{2^N},\frac{v2^{B^*}}{2^N})=
   < M-G| \delta_v \otimes \gamma_u >  \\-
   (\frac{u2^{A^*}}{2^N} \times \frac{v2^{B^*}}{2^N} ) 
   < M-G| e^N_0 \otimes e^N_0 > 
\end{eqnarray*}
where $M$ is defined by (\ref{disemp}) and where $G$ is the 
$\R^{2^N} \otimes \R^{2^N}$ vector defined by 
\begin{eqnarray*}
G=(Z^{0,0}_{0,0},Z^{0,0}_{0,1},\ldots,Z^{0,0}_{0,2^N-1},
   Z^{0,1}_{0,0},Z^{0,1}_{0,1},\ldots,Z^{0,1}_{0,2^N-1}, \cdots, \\
   Z^{0,2^N-1}_{0,0},Z^{0,2^N-1}_{0,1},\ldots,Z^{0,2^N-1}_{0,2^N-1}). 
\end{eqnarray*}
We have to expand the vectors $\gamma_u$ et $\delta_v$ on an appropriate basis.
Let $e_k^j$ be
the $\R^{2^N}$ vector associated with \( \{ k2^j+1,\ldots,(k+1)2^j \} \)
($0\leq j \leq N, \; 0\leq k \leq 2^{N-j}-1$).
Set \( \tilde{e}_k^j = e_k^j-e_{k+1}^j \)
for \( j\in \{0,\ldots,N-1 \} \), \( k\in \{0,\ldots,2^{N-u}-1 \} \),
$k$ even.
Thus \( {\cal B}=(e_0^N, \tilde{e}_k^j; k\in \{0,\ldots,N-1 \};
k \in \{ 0,\ldots,2^{N-u}-1 \}; k \mbox{ even}) \) is an orthogonal
basis of $\R^{2^N}$ and we have
\[ \gamma_u= \sum_{j=A^*}^{N-1} c^j_u \tilde{e}_{k(j,u)}^j  +
    \frac{u2^{A^*}}{2^N} e_0^N \]
where $k(j,u)$ is the only even integer such that
\[ u2^{A^*} \in ]k(j,u)2^j,(k(j,u)+2)2^j] \]
and where
\[ c^j_u=\frac{<\gamma_u|\tilde{e}_{k(j,u)}^j>}{2^{j+1}}. \]
In the same way, we have
\[  \delta_v = \sum_{i=B^*}^{N-1} c^i_l \tilde{e}_{l(i,v)}^i  +
    \frac{v2^{B^*}}{2^N} e_0^N
\]
where $l(i,v)$ is the only even integer such that
\[ v2^{B^*} \in ](l(i,v)-1)2^i,(l(i,v)+1)2^i] \]
and where
\[ c^i_v=\frac{<\delta_v|\tilde{e}_{l(i,v)}^i>}{2^{i+1}}. \]
The properties of coefficients $c^i_v,c^j_u$ will be useful throughout
this paper, therefore we state these properties by Lemma \ref{coeff}.
\begin{lem} \label{coeff} a) $0\leq c^i_v,c^j_u \leq 1/2$, \\
b) if  $i\geq B$ we have $c^i_v \leq \displaystyle \frac{2^B}{2^{i+1}}$, \\
c) if  $j\geq A$ we have $c^j_u \leq \displaystyle \frac{2^A}{2^{j+1}}$.
\end{lem}
Using the previous expansions we obtain
\begin{eqnarray*} 
   < M-G| \delta_v \otimes \gamma_u > -
   (\frac{u2^{A^*}}{2^N} \times \frac{v2^{B^*}}{2^N} ) 
   < M-G| e^N_0 \otimes e^N_0 > =
   \sum_{i=B^*}^{N-1} \sum_{j=A^*}^{N-1} c^i_v c^j_u 
  < M-G|\tilde{e}_{l(i,v)}^i \otimes \tilde{e}_{k(j,u)}^j> \\
\\ 
+\frac{v2^{B^*}}{2^N}< M-G| e^N_0 \otimes (\gamma_u-\frac{u2^{A^*}}{2^N}e^N_0)>
+\frac{u2^{A^*}}{2^N}< M-G| (\delta_v-\frac{v2^{B^*}}{2^N}e^N_0)\otimes e^N_0>.
\end{eqnarray*}
Let us recall that
\( C_1= \displaystyle { \frac{3 \T }{2} + 2 \left( \frac{3 C_0 }{4} \right)^{2/3} } \) 
when $ab=1$  and
\( C_1= \displaystyle {  \frac{3 \T }{2} +  2 \left( 3 C_0 \right)^{2/3} } \)
when $ab<1$.      
Let $Q$ be the probability to be controlled to obtain
Lemma \ref{principal} :
\begin{eqnarray*} Q=\mathop{\mathbb{P}} \left( 
  \max_{1\leq u \leq 2^{A-A^*}} 
  \max_{1\leq v \leq 2^{B-B^*}}    
  |nH_n(\frac{u2^{A^*}}{2^N},\frac{v2^{B^*}}{2^N})-
  n\frac{u2^{A^*}}{2^N}\frac{v2^{B^*}}{2^N}
  -\sqrt{n} D(\frac{u2^{A^*}}{2^N},\frac{v2^{B^*}}{2^N})| \right. \\ \geq  
  \left. {0.6 (x+C_1\log (nab))^{3/2}}
  \right) . 
\end{eqnarray*}
We have :
\begin{eqnarray*}
Q \leq \proba{ \max_{1\leq u \leq 2^{A-A^*}} 
  \max_{1\leq v \leq 2^{B-B^*}}   
  |\sum_{i=B^*}^{N-1} \sum_{j=A^*}^{N-1} c^i_v c^j_u 
  < M-G|\tilde{e}_{l(i,v)}^i \otimes \tilde{e}_{k(j,u)}^j> |
   \geq 0.6 ( 0.8 x +1.5 \T \log (nab))^{3/2} } \\
   + \proba{ \frac{2^B}{2^N} \max_{1\leq u \leq 2^{A-A^*}}
   |nH_n(\frac{u2^{A^*}}{2^N},1)-n\frac{u2^{A^*}}{2^N}-
   \sqrt{n} D(\frac{u2^{A^*}}{2^N},1)| \geq 
   0.6 (0.1 x + (\theta C_0)^{2/3} \log (nab))^{3/2} } \\
   + \proba{ \frac{2^A}{2^N}  \max_{1\leq v \leq 2^{B-B^*}}
   |nH_n(1,\frac{v2^{B^*}}{2^N})-n\frac{v2^{B^*}}{2^N}-
   \sqrt{n} D(1,\frac{v2^{B^*}}{2^N})| \geq 
   0.6 (0.1 x + (\theta C_0)^{2/3} \log (nab))^{3/2} }
\end{eqnarray*}
with $\theta=3$ when $ab<1$ and $\theta=3/4$ when $ab=1$.
The two last terms are completely analogous.
We detail the upper bound for the last term in the case $ab<1$. 
In this case, we use Inequality (\ref{in1ter}) and the relations 
$2^{A-N} < 2a $, $2^{B-N}<2b$, 
\( (\log (nab))/(2a) \geq (\log (2nb))/4 \),
and we obtain, considering $C_0 \geq 12$,
\begin{eqnarray*}
\proba{ \frac{2^A}{2^N}  \max_{1\leq v \leq 2^{B-B^*}}
   |nH_n(1,\frac{v2^{B^*}}{2^N})-n\frac{v2^{B^*}}{2^N}-
   \sqrt{n} D(1,\frac{v2^{B^*}}{2^N}) | \geq 
  0.6 (0.1+ (3C_0)^{2/3} \log (nab))^{3/2} } \\
   \leq \proba{ \sup_{0 \leq s \leq 2b} 
   |nG_n(s)-ns-\sqrt{n} D(1,s)| \geq 0.31 x + C_0 \log (2nb) } \\
   \leq \Lambda_0 \exp (-\lambda_0 x /5) .
\end{eqnarray*}
Considering moreover the relation
\[ 2^{A-A^*} 2^{B-B^*} \leq \frac{2 (nab)^2}{31 \T \log (nab) }, \]
the proof of Lemma \ref{principal} is achieved with the following lemma :
\begin{lem} \label{levoila}
In the case $nab > 496$, for all
\( u \in \{1, \ldots, 2^{A-A^*} \} \) and for all
\( v \in \{1, \ldots, 2^{B-B^*} \} \) we have :
\begin{eqnarray*}
   \proba {    |\sum_{i=B^*}^{N-1} \sum_{j=A^*}^{N-1} c^i_v c^j_u 
   < M-G|\tilde{e}_{l(i,v)}^i \otimes \tilde{e}_{k(j,u)}^j> |
   \geq ((x/2)+\T \log (nab))^{3/2} }  \\ \leq \Lambda_4 \log (nab) 
   \exp (-\lambda_4 x -2 \log (nab) )
\end{eqnarray*}
where $\Lambda_4 ,\lambda_4$ are absolute positive constants.
\end{lem}
Let $T(u,v)$ be the term to be controlled :
\begin{equation} T(u,v)=\sum_{i=B^*}^{N-1} \sum_{j=A^*}^{N-1} c^i_v c^j_u 
   < M-G|\tilde{e}_{l(i,v)}^i \otimes \tilde{e}_{k(j,u)}^j>.
\label{term} \end{equation}
The control of $T(u,v)$ is obtained from an exponential inequality of
Van de Geer (1995) and de la Pe\~{n}a (1999). 
This inequality, to which we devote Section \ref{rajout},
allows to control some martingales on an appropriate event.
The control of $T(u,v)$ will be of type
\begin{eqnarray*}  \proba { |T(u,v)| \geq (x/2+\T \log (nab))^{3/2} } \\
\leq \proba { \{ |T(u,v)| \geq (x/2+\T \log (nab))^{3/2} \} \cap \Theta(u,v) }
+ \proba { (\Theta(u,v))^{\mbox{c}} }. 
\end{eqnarray*}
We define below the event $\Theta(u,v)$. \\
For technical reasons, we have to consider some events where 
$U^{i,l}_{j,k}$ is close of $\esp { U^{i,l}_{j,k} }=\gamma 2^{i+j-N}$.
These events are of type
\begin{equation} 
{\cal E}_{j,k}^{i,l} = \{ | U^{i,l}_{j,k}-\gamma 2^{i+j-N} | \leq 
    \epsilon \gamma 2^{i+j-N}  \}.  \label{small}
\end{equation}
We take from now on  $\epsilon=1/2$ . 
The events $({\cal E}_{j,k}^{i,l})^{\mbox{c}}$ are controled in probability
by the following lemma (Benett (1962) and Wellner(1978), see also
Cs\"{o}rg\H{o} et Horv\'ath (1993) page 116) :
\begin{lem} \label{controle}
Let $Z$ be a binomial variable with expectation $m$. Then, for any positive $y$ 
and for any sign $\epsilon$ we have 
\( \proba { \epsilon (Z-m) \geq y } \leq \exp( -m h(y/m)) \) where the function
$h$ is defined for $t>-1$ by \( h(t)=(1+t) \ln(1+t) -t \).
\end{lem}
Thus we obtain
\[ \proba{ ({\cal E}_{j,k}^{i,l})^{\mbox{c}} } \leq 2 
   \exp (-\gamma 2^{i+j-N} h(\epsilon) )
\]
and we see that this control is suitable only when
$2^{i+j-N}$ is of order ${x}+C \log (nab)$. Therefore
we define the integers $M(i)$ and ${\cal M}(j)$ by :
\[ M(i)=\left\{ \begin{array}{ll} B^*+A-i-2 
        & \mbox{ for }i=B^*,\ldots,B-1 \\
        A^*-1 & \mbox{ for }i \geq B-1.
        \end{array} \right.
\]
\[ {\cal M}(j)=\left\{ \begin{array}{ll} A^*+B-j-2 
       & \mbox{ for }j=A^*,\ldots,A-1 \\
       B^*-1 & \mbox{ for }j \geq A-1.
       \end{array} \right.
\]
We have \( A^*-1 \leq M(i) \leq A-2 \), 
\( B^*-1 \leq {\cal M}(j) \leq B -2 \) and  
\[ \frac{ (x/2)+\T \log(nab) }{2\gamma }
    \leq 2^{i+M(i)-N}= 2^{j+{\cal M}(j)-N} <
\frac{ (x/2)+\T \log(nab) }{\gamma }. \]
We define the event $\Theta_0(u,v)$ by :
\begin{equation} \label{evt0}
\Theta_0(u,v)=\bigcap_{i=B^*}^{N-1} \bigcap_{j=M(i)+1}^{N-1}
\left( {\cal E}_{j,k(j,u)}^{i+1,l(i,v)/2} \cap
{\cal E}_{j,k(j,u)+1}^{i+1,l(i,v)/2} \cap
{\cal E}_{j+1,k(j,u)/2}^{i,l(i,v)} \cap
{\cal E}_{j+1,k(j,u)/2}^{i,l(i,v)+1}  \right) \end{equation}
where the basic event
${\cal E}_{j,k}^{i,l}$ is defined by (\ref{small}).
We define the event $\Theta_1(u,v)$ by :
\begin{eqnarray} \label{evt1}
\Theta_1(u,v)=\bigcap_{i=B^*}^{N-1} \left\{ \sum_{j=M(i)+1}^{N-1}
(\alpha_j \beta_i)^{1/2} \frac{ \left( {U}_{j,k(j,u)}^{i+1,l(i,v)/2}
-{U}_{j,k(j,u)+1}^{i+1,l(i,v)/2} \right)^2 }{
{U}_{j+1,k(j,u)/2}^{i+1,l(i,v)/2} } \leq (x/2) + \T \log (nab) \right\} \\
\cap 
\bigcap_{j=A^*}^{N-1} \left\{ \sum_{i={\cal M}(j)+1}^{N-1}
(\alpha_j \beta_i)^{1/2} \frac{ \left( {U}_{j+1,k(j,u)/2}^{i,l(i,v)}
-{U}_{j,+1k(j,u)/2}^{i,l(i,v)+1} \right)^2 }{
{U}_{j+1,k(j,u)/2}^{i+1,l(i,v)/2} } \leq  (x/2) + \T \log (nab) \right\}
\nonumber
\end{eqnarray}
with \[ \alpha_j=\inf (\frac{1}{2},\frac{2^A}{2^{j+1}}) , \]
and in the same way, 
\[ \beta_i=\inf (\frac{1}{2},\frac{2^B}{2^{i+1}}) . \]
The event $\Theta(u,v)$ on which we can control $T(u,v)$ is defined by
\begin{equation} \label{evt} \Theta(u,v)=\Theta_0(u,v) \cap \Theta_1(u,v). 
\end{equation}
Thus the proof of Lemma \ref{levoila},
and consequently the proof of Lemma \ref{principal},
is achieved with the two following lemmas :
\begin{lem} \label{contevt} Let $\Theta(u,v)$ be the event defined by
(\ref{evt}). In the case $nab>496$,
for all \( u \in \{1, \ldots, 2^{A-A^*} \} \) and for all
\( v \in \{1, \ldots, 2^{B-B^*} \} \) we have
\[ \proba { ( \Theta (u,v) )^{\mbox{c}} } \leq 
    \Lambda_5  \log(nab) \exp (-\lambda_5 x-2 \log (nab) ) 
\]
where $\Lambda_5 ,\lambda_5$ are absolute positive constants.
\end{lem}
\begin{lem} \label{conterm} Let $\Theta(u,v)$ be the event defined by
(\ref{evt}) and let $T(u,v)$ be the term defined by (\ref{term}). In the case 
$nab>496$, for all \( u \in \{1, \ldots, 2^{A-A^*} \} \) and for all 
\( v \in \{1, \ldots, 2^{B-B^*} \} \) we have
\[  \proba { \{ |T(u,v)| \geq ((x/2)+\T \log (nab))^{3/2} \} 
    \cap \Theta(u,v) }
    \leq \Lambda_6  \exp (-\lambda_6 x- 2\log(nab) ) 
\]
where $\Lambda_6 ,\lambda_6$ are absolute positive constants.
\end{lem}
The term $T(u,v)$ may be written as \[ T(u,v) = T_1(u,v)+ T_2(u,v) \]
where the terms $T_1(u,v),T_2(u,v)$ are defined by
\begin{eqnarray} \label{term1}
T_1(u,v)=\sum_{i=B^*}^{B-2} \sum_{j=A^*}^{M(i)} c^i_v c^j_u 
   < M-G|\tilde{e}_{l(i,v)}^i \otimes \tilde{e}_{k(j,u)}^j>, \\
\label{term2}
T_2(u,v)=\sum_{i=B^*}^{N-1} \sum_{j=M(i)+1}^{N-1} c^i_v c^j_u 
   < M-G|\tilde{e}_{l(i,v)}^i \otimes \tilde{e}_{k(j,u)}^j>.
\end{eqnarray}
Thus Lemma \ref{conterm} is proved by the two following lemmas :
\begin{lem} \label{conterm1} Let $\Theta_0(u,v)$ be the event defined by
(\ref{evt0}) and let $T_1(u,v)$ be the term defined by (\ref{term1}). 
In the case 
$nab>496$, for all \( u \in \{1, \ldots, 2^{A-A^*} \} \) and for all 
\( v \in \{1, \ldots, 2^{B-B^*} \} \) we have
\[  \proba { \left\{ |T_1(u,v)| \geq 
   \frac{ ( (x/2)+\T \log (nab) )^{3/2} }{2}  \right\} 
\cap \Theta_0(u,v) }
    \leq \Lambda_7  \exp (-\lambda_7 x- 2\log(nab) ) 
\]
where $\Lambda_7 ,\lambda_7$ are absolute positive constants.
\end{lem}
\begin{lem} \label{conterm2} Let $\Theta(u,v)$  be the event defined by (\ref{evt}), 
and let $T_2(u,v)$  be the term defined by (\ref{term2}). 
In the case 
$nab>496$, for all \( u \in \{1, \ldots, 2^{A-A^*} \} \) and for all 
\( v \in \{1, \ldots, 2^{B-B^*} \} \) we have
\[  \proba { \left\{ |T_2(u,v)| \geq \frac{( (x/2)+\T \log (nab))^{3/2} }{2}
\right\} 
\cap \Theta(u,v) }
    \leq \Lambda_8  \exp (-\lambda_8 x- 2\log(nab) ) 
\]
where $\Lambda_8 ,\lambda_8$ are absolute positive constants.
\end{lem}
{\bf Conclusion:} The proof of Lemma \ref{levoila}, and consequentely the proof of
Lemma \ref{principal}, is achieved with Lemmas
\ref{contevt}, \ref{conterm1}, \ref{conterm2}. 
We prove Lemma \ref{contevt} by Section \ref{5}, 
Lemma \ref{conterm1} by Section \ref{termosc}, 
Lemma \ref{conterm2} by Section \ref{der}, Section \ref{rajout} is
devoted to the result of Van de Geer (1995) and de la Pe\~na (1999).
\section{Proof of Lemma \ref{contevt}.}
\label{5}
By Lemma \ref{controle} we obtain
\begin{eqnarray*} \proba { \left( \Theta_0(u,v) \right)^{\mbox{c}} }
\leq 8 \sum_{i=B^*}^{B-2}  \sum_{j=M(i)+1}^{N-1} \exp \left(
-\gamma 2^{i+j+1-N} h(\epsilon )  \right) 
+8\sum_{i=B-1}^{N-1}  \sum_{j=A^*}^{N-1}
\exp \left(-\gamma 2^{i+j+1-N} h(\epsilon ) \right) 
\\
\leq \frac{ 8 \log (nab)}{\log (2) }
\sum_{s \geq 0} \exp \left(
-{2} h(\epsilon ) ((x/2)+\T \log(nab))  2^s \right) 
+ 8 \sum_{r \geq 0}  \sum_{s \geq 0} \exp \left(
-{2} h(\epsilon )((x/2)+\T \log(nab))  2^r 2^s \right) 
\\
\leq R_0 \exp(-\gamma_0 x+2 \log(nab) ) 
\end{eqnarray*} 
where $R_0,\gamma_0$ are absolute positive constants.
We set 
\[  \Delta^{i+1,l(i,v)/2}_{j,k(j,u)}= (\alpha_j \beta_i )^{1/2}
                   \frac{ \left( U^{i+1,l(i,v)/2}_{j,k(j,u)}
                   -U^{i+1,l(i,v)/2}_{j,k(j,u)+1} \right)^2 }
                   {U^{i+1,l(i,v)/2}_{j+1,k(j,u)/2}}                   
\]
\[                  
      \tilde{\Delta}^{i,l(i,v)}_{j+1,k(j,u)/2}=(\alpha_j \beta_i )^{1/2}
      \frac{ \left( U^{i,l(i,v)}_{j+1,k(j,u)/2}-U^{i,l(i,v)+1}
      _{j+1,k(j,u)/2} \right)^2 }{U^{i+1,l(i,v)/2}_{j+1,k(j,u)/2}}.
\]
With these notations, we have
\begin{eqnarray*}
\Theta_1(u,v)=\bigcap_{i=B^*}^{N-1} \left\{ \sum_{j=M(i)+1}^{N-1}
 \left| \Delta^{i+1,l(i,v)/2}_{j,k(j,u)} \right| \leq (x/2) + \T \log (nab) \right\} \\
\cap 
\bigcap_{j=A^*}^{N-1} \left\{ \sum_{i={\cal M}(j)+1}^{N-1}
\left| \tilde{\Delta}^{i,l(i,v)}_{j+1,k(j,u)/2} \right| \leq  (x/2) + \T \log (nab) \right\}
\end{eqnarray*}
and
\begin{eqnarray*} \proba { \left( \Theta_1(u,v) \right)^{\mbox{c}} 
\cap \Theta_0(u,v) } \leq 
\sum_{i=B^*}^{N-1} \proba { \left\{ \sum_{j=M(i)+1}^{N-1}
\Delta^{i+1,l(i,v)/2}_{j,k(j,u)} > (x/2) + \T \log (nab) \right\} 
\cap  \Theta_0(u,v) }\\+ \sum_{j=A^*}^{N-1} 
\proba { \left\{ \sum_{i={\cal M}(j)+1}^{N-1}
\tilde{\Delta}^{i,l(i,v)}_{j+1,k(j,u)/2} >  (x/2) + \T \log (nab) \right
\}\cap \Theta_0(u,v) }
\end{eqnarray*}
We use the first inequality of Tusn\'ady's Lemma (1977 b)
(conditional construction of a multinomial vector) 
proved in 1989 by Bretagnolle and Massart.
In order to apply this lemma directly, we express it with our notations.
\begin{lem} \label{tusna} (Tusn\'ady)
For all \( i \in \{ B^*, \ldots, N-1 \} \)
there exists i.i.d. ${\cal N}(0,1)$ random variables, denoted by 
\( {\xi}^{i+1,l(i,v)/2}_{j,k(j,u)};j=A^*,\ldots,  N-1 \), such that
\[ \left|  U^{i+1,l(i,v)/2}_{j,k(j,u)}
                   -U^{i+1,l(i,v)/2}_{j,k(j,u)+1} \right|
    \leq 2 \left( 1 +\frac{ \sqrt{ U^{i+1,l(i,v)/2}_{j+1,k(j,u)/2} } }{2}  
    \left| {\xi}^{i+1,l(i,v)/2}_{j,k(j,u)} \right| \right) .
\]
In the same way, for all \( j \in \{ A^*, \ldots, N-1 \} \)
there exists i.i.d. ${\cal N}(0,1)$ random variables, denoted by 
\( \tilde{{\xi}}^{i,l(i,v)}_{j+1,k(j,u)/2};i=B^*,\ldots,  N-1 \), such that
\[ \left|  U^{i,l(i,v)}_{j+1,k(j,u)/2}
                   -U^{i,l(i,v)+1}_{j+1,k(j,u)/2} \right|
    \leq 2 \left( 1 +\frac{ \sqrt{ U^{i+1,l(i,v)/2}_{j+1,k(j,u)/2} } }{2}  
    \left| \tilde{\xi}^{i,l(i,v)}_{j+1,k(j,u)/2} \right| \right) .
\]
\end{lem}
Lemma \ref{tusna} yields that on $\Theta_0(u,v)$ we have :
\begin{eqnarray*}
\left| {\Delta}^{i+1,l(i,v)/2}_{j,k(j,u)} \right| \leq 8 \sqrt{\beta_i}
\sum_{j=M(i)+1}^{N-1} \sqrt{\alpha_j} \left( 
\frac{1}{{U}^{i+1,l(i,v)/2}_{j+1,k(j,u)/2}}+0.25 
\left( {\xi}^{i+1,l(i,v)/2}_{j,k(j,u)} \right)^2 \right) \\
\leq \frac{4}{\T \log (nab)} + \sqrt{\beta_i}
\sum_{j=M(i)+1}^{N-1} 0.5 \sqrt{\alpha_j}
\left( {\xi}^{i+1,l(i,v)/2}_{j,k(j,u)} \right)^2 ,
\end{eqnarray*}
and also
\[ \left| \tilde{\Delta}^{i+1,l(i,v)/2}_{j,k(j,u)} \right|
\leq \frac{4}{\T \log (nab)} + \sqrt{\alpha_j}
\sum_{i={\cal M}(j)+1}^{N-1} 0.5 \sqrt{\beta_i} 
\left( \tilde{\xi}^{i,l(i,v)}_{j+1,k(j,u)/2} \right)^2 .
\]
Hence, since $\T \geq 10$ and  $nab>496$, by setting
$\TT=9.9$, one gets :
\begin{eqnarray*} \proba { \left( \Theta_1(u,v) \right)^{\mbox{c}} 
\cap \Theta_0(u,v) } \leq 
\sum_{i=B^*}^{N-1} \proba {  \sum_{j=M(i)+1}^{N-1}
\sqrt{\alpha_j} \left( {\xi}^{i+1,l(i,v)}_{j,k(j,u)} \right)^2
> \frac{ 2 ((x/2) + \TT \log (nab))}{ \sqrt{\beta_i} } }
\\+ \sum_{j=A^*}^{N-1} 
\proba {  \sum_{i={\cal M}(j)+1}^{N-1}
\sqrt{\beta_i} \left( {\xi}^{i,l(i,v)}_{j+1,k(j,u)/2} \right)^2 
 > \frac{ 2 ((x/2) + \TT \log (nab))}{ \sqrt{\alpha_j} } }.
\end{eqnarray*}
The control of the two terms being completely analogous, we obtain
\[ \proba { \left( \Theta_1(u,v) \right)^{\mbox{c}} 
\cap \Theta_0(u,v) } \leq 2
\sum_{i=B^*}^{N-1} \proba {  \sum_{j=M(i)+1}^{N-1}
\sqrt{\alpha_j} \left( {\xi}^{i+1,l(i,v)}_{j,k(j,u)} \right)^2
> \frac{ 2 ((x/2) + \TT \log (nab))}{ \sqrt{\beta_i} } }.
\]
We use Cramer-Chernov Inequality :
\begin{lem} \label{chi}
Let \( \zeta_1,\ldots,\zeta_d \) be i.i.d. 
${\cal N}(0,1)$ random variables and let
\( \lambda_1,\ldots,\lambda_d \) be some positive integers. We have 
\[ \proba{ \sum_{i=1}^d \lambda_i \zeta_i^2 \geq z }
    \leq \inf_{0<r<\inf_i 1/(2\lambda_i)} \exp (-rz -\frac{1}{2}
    \sum_{i=1}^d\ln(1-2\lambda_i r)) . \]
\end{lem}
We take $r=1/2$ and we use \( \ln (1-x) \geq -1.8 x \) for
$x \leq 1/\sqrt{2}$. We get
\[ \proba { \left( \Theta_1(u,v) \right)^{\mbox{c}} 
\cap \Theta_0(u,v) } \leq 2
\sum_{i=B^*}^{N-1} \exp \left(
- \frac{ 2 ((x/2) + \TT \log (nab))}{ \sqrt{\beta_i} } +0.9 (A-M(i)+1) 
+\frac{ \sqrt{2}+1}{\sqrt{2}-1} \right) .
\]
We conclude with \( A-M(i)+1 \leq A^*-A+2 \leq (\log (nab))/(\log(2)) \) :
\begin{eqnarray*} \proba { \left( \Theta_1(u,v) \right)^{\mbox{c}} 
\cap \Theta_0(u,v) } \leq R_1 (nab)^{1.3} [ \;
(B-B^*+1) \exp \left(-\sqrt{2} ((x/2) + \TT \log (nab)) \right) 
\\  +\sum_{s\geq 0} \exp \left( -\sqrt{2^s} ((x/2) + \TT \log (nab)) 
\right) \; ]
\\
\leq R_1 \log (nab) \exp ( -\gamma_1 x- 2 \log (nab) )
\end{eqnarray*}
where $R_1,\gamma_1$ are absolute positive constants.
\section{Exponential inequality for martingales.} \label{rajout}
We are devoting a section to this inequality, because we use it greatly 
throughout the proof of Lemmas  \ref{conterm1} and \ref{conterm2} 
(Sections \ref{termosc} and \ref{der}). 
Van de Geer in 1995, then de la Pe\~na in 1999, have generalized Bernstein
Inequality to some not bounded martingales. 
It turned out (and this is rather surprising) that the error terms emanating
from Hungarian constructions
(in this paper, this is the term $T(u,v)$ in Lemma \ref{conterm}) 
are not bounded martingales exactly verifying  assumptions
of Van de Geer's or de la Pe\~na's Theorem. 
All Hungarian constructions of a dimension larger than 1
may probably be dealt with from this new point of view.
In this paper, we use de la Pen\~{a}'s notations.
First we recall his theorem, then we express it in a 
form appropriate to this paper.
\begin{theo} \label{pena} (Van de Geer, de la Pe\~{n}a) Let $(d_j)$ be a sequence 
adapted to the increasing filtration $(F_j)$ with $\esp { d_j/F_{j-1} } =0$,
$\esp { (d_j)^2/F_{j-1} } =\sigma_j^2$, 
${\cal V}_T^2=\sum_{j=1}^T \sigma_j^2$. Assume that
\begin{equation} \label{cond1}
\esp { |d_j|^k/F_{j-1} } \leq \frac{k!}{2} c^{k-2} \sigma_j^2 
\mbox{ \hspace{ 1cm} p.s.} 
\end{equation}
or 
\[ \proba { |d_j| \leq c }=1 \]
for $k>2$, $0<c<\infty $. Then, for all $x,y>0$,
\[ \proba { \sum_{j=1}^T d_j \geq x, \;{\cal V}_T^2 \leq y \mbox{ for 
some } T } \leq \exp(-\frac{x^2}{2(y+cx)}). \]
\end{theo} 
\begin{lem} \label{simp}
Let $(d_j)$ and $(F_j)$ be defined by Theorem \ref{pena}. If
\begin{equation} \label{cond2} 
\esp { |d_j|^k/F_{j-1} } \leq \frac{k!}{2} c^k 
\end{equation}
for $k\geq 2$, $0 < c < \infty$, then the condition (\ref{cond1}) holds.
\end{lem}
\underline{Proof of Lemma \ref{simp}} :
\begin{eqnarray*} \esp { |d_j|^k/F_{j-1} }  = 
\esp { \left( |d_j|^k \ind \left\{ |d_j|^2 \leq c^2 \right\} \right) /F_{j-1} }
+\esp { \left( |d_j|^k \ind \left\{ |d_j|^2 > c^2 \right\} \right) /F_{j-1} } 
\\
{ \leq   c^{k-2} \esp { \left( |d_j|^2 
\ind \left\{  \esp { |d_j|^2 / F_{j-1} } \leq c^2 \right\} \right) /F_{j-1} }
+\esp { \left( |d_j|^k 
\ind \left\{ \esp { |d_j|^2 /F_{j-1} } > c^2 \right\} \right) /F_{j-1} } }\\
{ \leq    c^{k-2} \esp {  |d_j|^2  /F_{j-1} }
\ind \left\{ \esp { |d_j|^2 / F_{j-1} } \leq c^2 \right\} 
+\frac{k!}{2} c^k \ind \left\{ \esp { |d_j|^2 /F_{j-1} } > c^2 \right\} }\\
{  \leq  c^{k-2} \esp {  |d_j|^2  /F_{j-1} }
\ind \left\{ \esp { |d_j|^2 / F_{j-1} } \leq c^2 \right\} 
+\frac{k!}{2} c^{k-2} \esp { |d_j|^2 /F_{j-1} }
\ind \left\{ \esp { |d_j|^2 /F_{j-1} } > c^2 \right\} }\\
\leq   \frac{k!}{2} c^{k-2} \esp { |d_j|^2 /F_{j-1} }.\\
\Box 
\end{eqnarray*}
Lemma \ref{simp} combined with Cauchy-Schwarz Inequality
gives the following Lemma :
\begin{lem} \label{simpli}
Let $(d_j)$ and  $(F_j)$ be defined by Theorem \ref{pena}. If 
\begin{equation} \label{cond3}
\esp { |d_j|^{2k}/F_{j-1} } \leq \frac{(2k)!}{2^k k!} c^{2k} 
\end{equation}
for $k\geq 1$, $0 < c < \infty$, then the condition (\ref{cond1}) holds.
\end{lem}
In Sections \ref{termosc} and \ref{der} we use Theorem
\ref{pena} in the following form :
\begin{theo} \label{penafin}
Let $(d_j)$ and $(F_j)$ be defined by Theorem \ref{pena}. 
Let $\Theta$ be an event such that on $\Theta$
we have (\ref{cond2}) or (\ref{cond3})
and ${\cal V}_T^2 \leq y$ where ${\cal V}_T^2$
is defined by Theorem \ref{pena}. Then, for all $x>0$,
\[ \proba { \{ \sum_{j=1}^T d_j \geq x \} \cap \Theta } 
\leq \exp(-\frac{x^2}{2(y+cx)}). \]
\end{theo}
\underline{Proof of Theorem \ref{penafin}} :
Let us denote by $E_j$ the event
\( \{ \esp {  |d_j|^{2k}/F_{j-1} } \leq 
\displaystyle  \frac{(2k)!}{2^k k!} c^{2k} 
\mbox{ for all } k\geq 1 \} \) or the event 
\( \{ \esp {  |d_j|^{k}/F_{j-1} } \leq 
\displaystyle  \frac{k!}{2} c^{2k} 
\mbox{ for all } k\geq 2 \} \). We have :
\begin{eqnarray*} \{ \sum_{j=1}^T d_j \geq x \} \cap \Theta 
&=& \{  \sum_{j=1}^T d_j \ind \{ E_j \} + \sum_{j=1}^T d_j 
\ind \{ E_j^{\mbox{c}} \}  \geq x \} \cap \Theta \\
&=& \{  \sum_{j=1}^T d_j \ind \{ E_j \} \geq x \} \cap \Theta \\
&\subset & \{  \sum_{j=1}^T d_j \ind \{ E_j \} \geq x \}
\cap \{ {\cal V}_T^2  \leq y \}
\end{eqnarray*}
and in the same way,
\[ \{ -\sum_{j=1}^T d_j \geq x \} \cap \Theta 
\subset  \{  -\sum_{j=1}^T d_j \ind \{ E_j \} \geq x \}
\cap \{ {\cal V}_T^2  \leq y \} .
\]
Then we apply Lemmas \ref{simp} or \ref{simpli} 
and Theorem \ref{pena}
to $(D_j,F_j)$ with $D_j=d_j\ind \{E_j \}$. \hspace{2cm}  $\Box $
\section{Proof of Lemma \ref{conterm1}.} \label{termosc}
Let $P_1(u,v)$ be the probability to be controlled
to get Lemma \ref{conterm1} :
\[
P_1(u,v)=\proba { |
\left\{ T_1(u,v)| \geq \frac{1}{2}((x/2)+\T \log (nab))^{3/2} \right\}
\cap \Theta_0(u,v) }.
\]
We separate the Gaussian and empirical parts :
\[ T_1(u,v) \leq T_1^E(u,v) + T_1^G(u,v), \]
where the terms $T_1^E(u,v)$, $T_1^G(u,v)$ are defined by
\begin{eqnarray*}
T_1^E(u,v)=\sum_{i=B^*}^{B-2} \sum_{j=A^*}^{M(i)} c^i_v c^j_u 
   < M|\tilde{e}_{l(i,v)}^i \otimes \tilde{e}_{k(j,u)}^j>, \\
T_1^G(u,v)=\sum_{i=B^*}^{B-2} \sum_{j=A^*}^{M(i)} c^i_v c^j_u 
   < G|\tilde{e}_{l(i,v)}^i \otimes \tilde{e}_{k(j,u)}^j>. 
\end{eqnarray*}
We have \[ P_1(u,v) \leq P_1^E(u,v) + P_1^G(u,v), \]
where the probabilities $P_1^E(u,v)$, $P_1^G(u,v)$ are defined by
\[ P_1^E(u,v)=
   \proba { \left\{ |T_1^E(u,v)| \geq 
   \frac{\lambda}{2}((x/2)+\T \log (nab))^{3/2}  \right\}
   \cap \Theta_0(u,v) } \] 
\mbox{ and }  
\[ P_1^G(u,v)=
   \proba {\left\{  |T_1^G(u,v)| \geq 
   \frac{(1-\lambda)}{2}(x+\T \log (nab))^{3/2} \right\}
   \cap \Theta_0(u,v) }, \] 
with $\lambda =1/2$.
\subsection{Control of $ P_1^G(u,v)$.}
The control of $P_1^G(u,v)$ is directly, on observing that
with the notations of Section \ref{2} we have
\[ <G|\tilde{e}_{l}^i \otimes \tilde{e}_{k(j,u)}^j>= 4V^{i,l}_{j,k(j,u)}. \]
Consequently $T_1^E(u,v)$ is a Gaussian variable with expectation 0
and with variance equal to
\[ \sum_{i=B^*}^{B-2} \sum_{j=A^*}^{M(i)} (c^i_v c^j_u)^2 4 \gamma^2 2^{i+j-N} .
\] This variance is bounded  (\(0 \leq  c^i_v,c^j_u \leq 1/2 \)) by
\( \gamma^2 (B-B^*-1)2^{A^*+B-N-3} \), thus by $(\tau (x+\T \log (nab)))^2$, 
where $\tau$ is a positive constant verifying $\tau \leq ( \T \ln 4)^{-1/2}$.
Then using the well known inequality
\[ \proba { Y \geq t } \leq \frac{1}{t\sqrt{2\pi}}
\exp (-t^2/2) , \] where $Y$ denotes a standard Gaussian variable, we obtain
\[
   P_1^G(u,v) \leq \frac{8\tau}{\sqrt{2\pi \T}} 
   \exp (-\frac{(x/2)+\T \log (nab)}
   {32 \tau^2}) \leq R_2 \exp ( -\gamma_2 x -2 \log(nab)) 
\]
where $R_2 ,\gamma_2$ are absolute positive constants
(we use $\T \geq 10$, thus constants do not depend on $\T$).
\subsection{Control of $ P_1^E(u,v) $.} \label{62}
The control of $P_1^E(u,v)$ is more complicated because the variables
$< M|\tilde{e}_{l(i,v)}^i \otimes \tilde{e}_{k(j,u)}^j>$ are not
independent. 
Let us recall that $k(M(i),u)$ is the only even integer such that
\[ k(M(i),u)2^{M(i)} < u2^{A^*} \leq ( k(M(i),u) +2) 2^{M(i)} . \]
Let us denote by $\alpha_u$ the vector associated, 
according to Section \ref{4},
to \( \{ k(M(i),u)2^{M(i)}+1,\ldots , u2^{A^*} \} \) and let us denote
by $\beta_u$ the vector associated to 
\( \{ u2^{A^*}+1,\ldots , ( k(M(i),u)+2) 2^{M(i)} \} \). The expansion of
$\alpha_u$ on the basis ${\cal B}$ (defined by Section \ref{4}) is :
\[
\alpha_u = \sum_{j=A^*}^{M(i)} c^j_u \tilde{e}_{k(j,u)}^j  +
\sum_{j=M(i)+1}^{N-1} \epsilon^j_u \; 
\left( \frac{ u2^{A^*}-k(M(i),u)2^{M(i)} }{2^{j+1}} \right) 
\tilde{e}^j_{k(j,u)} 
+ \frac{ u2^{A^*}-k(M(i),u)2^{M(i)} }{2^{N}} e^N_0
\]
where $\epsilon^j_u$ is a sign defined by
\[   \epsilon^j_u=\left\{ \begin{array}{ll}
                  +1 & \mbox{ si } <\alpha_u|\tilde{e}^j_{k(j,u)}> > 0 \\
                  -1 & \mbox{ si } <\alpha_u|\tilde{e}^j_{k(j,u)}> < 0.
                  \end{array} \right.
\]
On the other hand, the expansion of $\alpha_u + \beta_u$ 
on the basis ${\cal B}$ is :
\[
\alpha_u + \beta_u = 
\sum_{j=M(i)+1}^{N-1} \epsilon^j_u \; 
\frac{ 2^{M(i)+1} }{2^{j+1}}
\tilde{e}^j_{k(j,u)} 
+ \frac{ 2^{M(i)+1} }{2^{N}} e^N_0.
\]
This gives
\[ \sum_{j=A^*}^{M(i)} c^j_u \tilde{e}_{k(j,u)}^j
    = \left( \frac{ (k(M(i),u)+2)2^{M(i)}-u2^{A^*} }{2^{M(i)+1}} \right) \alpha_u
    +\left( \frac{ u2^{A^*}-k(M(i),u)2^{M(i)} }{2^{M(i)+1}}\right)  \beta_u \]
and thus one obtains
\[ T_1^E(u,v) \leq d_1|\sum_{i=B^*}^{B-2} c^i_v 
   < M|\tilde{e}_{l(i,v)}^i \otimes \alpha_u >|+
   d_2|\sum_{i=B^*}^{B-2} c^i_v 
   < M|\tilde{e}_{l(i,v)}^i \otimes \beta_u >|
\]
with $d_1+d_2=1$. Finally, we have 
\begin{equation*}  P_1^E(u,v) \leq P_{1,\alpha}^E(u,v)
+ P_{1,\beta}^E(u,v), 
\end{equation*} 
where the probabilities $P_{1,\alpha}^E$, $ P_{1,\beta}^E$ are defined by
\begin{eqnarray*}
P_{1,\alpha}^E= \proba { \left\{ |\sum_{i=B^*}^{B-2} c^i_v 
   < M|\tilde{e}_{l(i,v)}^i \otimes \alpha_u >| \geq 
   \frac{(1-\lambda)}{2} 
   ((x/2)+\T \log (nab))^{3/2}  \right\} 
   \cap \Theta_0(u,v) }, \\
P_{1,\beta}^E= \proba { \left\{ |\sum_{i=B^*}^{B-2} c^i_v 
   < M|\tilde{e}_{l(i,v)}^i \otimes \beta_u >| \geq 
   \frac{(1-\lambda)}{2} ((x/2)+\T \log (nab))^{3/2}  \right\}
   \cap \Theta_0(u,v) }.
\end{eqnarray*}
We detail only the control of $P_{1,\alpha}^E(u,v)$ but 
the control of $P_{1,\beta}^E(u,v)$ is completely analogous.
First we verify the conditions of Theorem \ref{penafin}. The sequence
\[ \left( c^i_v < M|\tilde{e}_{l(i,v)}^i \otimes \alpha_u > \right),
i=B-2,\ldots,B^*
\]
is adapted to the decreasing filtration
\[ {\cal F}_0^{B-2} \subset {\cal F}_0^{B-3} \subset \ldots \subset
{\cal F}_0^{B^*}, \] because the variable
\( c^i_v < M|\tilde{e}_{l(i,v)}^i \otimes \alpha_u >  \) is
${\cal F}^{i}_0$ measurable
(${\cal F}_0^{i}$ is defined by Section \ref{2}).
Moreover, 
\begin{equation} \label{loi} 
{\cal L} \left( < M|{e}_{l(i,v)}^i \otimes \alpha_u > / {\cal F}^{i+1}_0 
\right) = {\cal B} \left( < M|{e}_{l(i,v)/2}^{i+1} \otimes \alpha_u > 
, \frac{1}{2} \right). 
\end{equation}
Let us recall that
\begin{equation} \label{bincent} < M|\tilde{e}_{l(i,v)}^i \otimes \alpha_u >=
   2 < M|{e}_{l(i,v)}^i \otimes \alpha_u > -
   < M|{e}_{l(i,v)/2}^{i+1} \otimes \alpha_u > .
\end{equation}
This yields 
\( \esp { < M|\tilde{e}_{l(i,v)}^i \otimes \alpha_u > / {\cal F}^{i+1}_0  }
=0 \).
As in Theorem \ref{pena}, let
\[ \left( \sigma^{i}(u,v)  \right)^2 = 
\esp { \left(  c^i_v < M|\tilde{e}_{l(i,v)}^i \otimes \alpha_u > 
\right)^2 / {\cal F}^{i+1}_0 } \mbox{ and }  
{\cal V}^2_{B^*}(u,v)=\sum_{i=B^*}^{B-2} \left( \sigma^{i}(u,v)  \right)^2 . \]
Using again (\ref{loi}) and  (\ref{bincent}) this gives
\[ \left( \sigma^{i}(u,v)  \right)^2 = (c^i_v )^2 
< M|{e}_{l(i,v)/2}^{i+1} \otimes \alpha_u > . \]
Since \( \displaystyle { \frac{k(M(i),u)}{2} \in \{ k(M(i)+1,u),
 k(M(i)+1,u)+1 \} } \)), on $\Theta_0(u,v)$ we have 
\begin{eqnarray}   < M|{e}_{l(i,v)/2}^{i+1} \otimes \alpha_u > 
\leq < M|{e}_{l(i,v)/2}^{i+1} \otimes \alpha_u + \beta_u > 
\leq \gamma (1+\epsilon) 2^{i+1+M(i)+1-N} =
\gamma (1+\epsilon) 2^{A^*+B-N}
\nonumber \\ \label{majvar} \leq 4(1+ \epsilon) ((x/2)+\T \log (nab)),
\end{eqnarray}
and this yields (using \( 0\leq c^i_v \leq 1/2 \)) 
that on $\Theta_0(u,v)$ we have :
\[
{\cal V}^2_{B^*}(u,v)
\leq (B-B^*-1) (1+\epsilon) ((x/2)+\T \log (nab)) 
\leq \frac{ (1+\epsilon)((x/2)+\T \log (nab))^2}{\T \ln (2)}.
\]
In order to verify condition (\ref{simpli}) we use the following lemma :
\begin{lem} \label{bineg} Let \( Z_1,\ldots ,Z_T \) be i.i.d. random variables,
\( \proba { Z_i=+1 }=\proba { Z_i=-1 }=1/2 \). 
We set $S=\sum_{i=1}^T Z_i$. For all $k\in \N^*$ we have
\[ \esp { S^{2k} } \leq \frac{ (2k)! }{ 2^k k! } T^k . \]
\end{lem}
\underline{Proof of Lemma \ref{bineg}} :
\[ S^{2k}=\sum_{i_1} \ldots \sum_{i_{2k}} Z_{i_1} \ldots Z_{i_{2k}}. \]
Let us define $N_w^{(i_1,\ldots,i_{2k})}$ as the number of indexes equal to $i_w$ :
\[ N_w^{(i_1,\ldots,i_{2k})}=\sum_{l=1}^{2k} \ind{ i_l=i_w }. \]
If there exists $w$ such that $N_w^{(i_1,\ldots,i_{2k})}$ is odd, then
\( \esp { Z_{i_1} \ldots Z_{i_{2k}} }=0 \). Thus
\begin{eqnarray*} \esp { S^{2k} } &=&
\sum_{ \{ (i_1,\ldots,i_{2k}) \mbox{ such that } N_w^{(i_1,\ldots,i_{2k})}
\mbox{ is even for all } w \in \{1,\ldots,2k \} \} }
\esp { Z_{i_1} \ldots Z_{i_{2k}} }
\\ & \\ &\leq & A \sum_{j_1} \ldots \sum_{j_{k}} 
\esp { Z_{j_1}^2 \ldots Z_{j_{k}}^2 },
\end{eqnarray*}
where \( A=\mbox{Card} \{ (i_1,\ldots,i_{2k}) \mbox{ such that } 
N_1^{(i_1,\ldots,i_{2k})}=\ldots =N_{2k}^{(i_1,\ldots,i_{2k})}=2 \} \) :
\[ A=\frac{ \mbox{C}^2_{2k} \mbox{C}^2_{2k-2} \ldots \mbox{C}^2_{2} }
     {k!}. \]
Since \( \esp { Z_{j_1}^2 \ldots Z_{j_{k}}^2 }=1 \) for all
$(j_1,\ldots,j_k)$, the proof is complete.
$\mbox{\hspace{3cm}} \Box$ \\ \\
Lemma \ref{bineg}, Equalities (\ref{loi}), (\ref{bincent}), 
the bound (\ref{majvar}) and the property \( 0 \leq c^i_v \leq 1/2 \)
yield (\ref{cond3}) :
\[
\esp { \left(  
c^i_v < M|\tilde{e}_{l(i,v)}^i \otimes \alpha_u > 
\right)^{2k} / {\cal F}^{i+1}_0 } \leq 
\frac{ (2k)! }{ 2^k k! } 
\left( \frac{< M|{e}_{l(i,v)/2}^{i+1} \otimes \alpha_u >}{4}  \right)^{k}
\leq \frac{ (2k)! }{ 2^k k! } c^{2k}
\]
with \[ c= \left( (1+\epsilon)  ((x/2)+\T \log (nab)) \right)^{1/2}. \]
We can now apply Theorem \ref{penafin} :
\begin{eqnarray*} P_{1,\alpha}^E (u,v) 
\leq 2 \exp 
\left( -
\frac{  ((x/2)+\T \log (nab))^{3} (1-\lambda)^2/4  }
{2 \left( 
\displaystyle
\frac{ (1+\epsilon)((x/2)+\T \log (nab))^2 }{ \T \ln (2)} +
\frac{\sqrt{(1+\epsilon)}(1-\lambda)}{2} ((x/2)+\T \log (nab))^{2} 
\right) 
}
\right) 
\\  \leq R_3 \exp (-\gamma_3 x+ 2 \log (nab) ) 
\end{eqnarray*}
where $R_3,\gamma_3$ are absolute positive constants
(we use $\T \geq 10$, thus constants do not depend on $\T$).
\section{Proof of Lemma \ref{conterm2}.} \label{der}
Let  $P_2(u,v)$ be the probability to be controlled 
to obtain Lemma \ref{conterm2} :
\[
P_2(u,v)=\proba { 
\left\{ |T_2(u,v)| \geq \frac{1}{2}((x/2)+\T \log (nab))^{3/2} \right\} 
\cap \Theta(u,v) }.
\]
Let us recall the definition of $T_2(u,v)$ :
\[
T_2(u,v)=\sum_{i=B^*}^{N-1} \sum_{j=M(i)+1}^{N-1} c^i_v c^j_u 
   < M-G|\tilde{e}_{l(i,v)}^i \otimes \tilde{e}_{k(j,u)}^j>.
\]
Lemma \ref{hyper} allows us to control on $\Theta(u,v)$
the expression
\begin{eqnarray*} \left| < M|\tilde{e}_{l(i,v)}^i \otimes \tilde{e}_{k(j,u)}^j>
      -\esp { < M|\tilde{e}_{l(i,v)}^i \otimes \tilde{e}_{k(j,u)}^j>
              / {\cal F}^i_{j+1} } \right. \\ \left. 
      - \sqrt{ \var { < M|\tilde{e}_{l(i,v)}^i \otimes \tilde{e}_{k(j,u)}^j>
              / {\cal F}^i_{j+1} } }
      \left( \frac{ \gamma 2^{i+j-N} }{4} \right)^{-1/2} 
      < G|\tilde{e}_{l(i,v)}^i \otimes \tilde{e}_{k(j,u)}^j>
   \right| 
\end{eqnarray*}
but unfortunately
$  \esp { < M|\tilde{e}_{l(i,v)}^i \otimes \tilde{e}_{k(j,u)}^j>
/ {\cal F}^i_{j+1} }  \neq 0 $. 
This is a flaw of the bivariate construction. Because of this flaw, in
dimension 2 the controls are more complicated than in dimension 1.
Perhaps another construction is conceivable, leading to the same theorem, 
but simpler and closer to the univariate construction.
This other construction is not still available, therefore we must write
the term $T_2(u,v)$ as a sum of three terms (instead of a sum of two
terms, which would be more natural). 
\\
Let us recall the notations of Lemma \ref{hyper} :  
\[  \delta^{i+1,l(i,v)/2}_{j,k(j,u)}=
                   \frac{U^{i+1,l(i,v)/2}_{j,k(j,u)}
                   -U^{i+1,l(i,v)/2}_{j,k(j,u)+1}}{U^{i+1,l(i,v)/2}_{j+1,k(j,u)/2}}                   
                   \;\;\; \mbox{ and } \;\;\; 
    \tilde{\delta}^{i,l(i,v)}_{j+1,k(j,u)/2}=
                   \frac{U^{i,l(i,v)}_{j+1,k(j,u)/2}-U^{i,l(i,v)+1}_{j+1,k(j,u)/2}}{U^{i+1,l(i,v)/2}_{j+1,k(j,u)/2}}.
\]
With the notations of Section \ref{2}, we have 
\[  < M|\tilde{e}_{l(i,v)}^i \otimes \tilde{e}_{k(j,u)}^j>
=U^{i,l(i,v)}_{j,k(j,u)}-U^{i,l(i,v)+1}_{j,k(j,u)}
-U^{i,l(i,v)}_{j,k(j,u)+1}+U^{i,l(i,v)+1}_{j,k(j,u)+1}
\]
and \[ < G|\tilde{e}_{l(i,v)}^i \otimes \tilde{e}_{k(j,u)}^j>
=4 V^{i,l(i,v)}_{j,k(j,u)}. \]
Hence, ones gets the expression
\begin{eqnarray*}
  < M-G |\tilde{e}_{l(i,v)}^i \otimes \tilde{e}_{k(j,u)}^j>= \\
  4 \left( U^{i,l(i,v)}_{j,k(j,u)}-
  \esp { U^{i,l(i,v)}_{j,k(j,u)} / {\cal F}^i_{j+1} } 
  - \sqrt{ \var { U^{i,l(i,v)}_{j,k(j,u)}  / {\cal F}^i_{j+1} } }               
      \left( \frac{ \gamma 2^{i+j-N} }{4} \right)^{-1/2} 
    V^{i,l(i,v)}_{j,k(j,u)} \right) \\
  +4 \left(  \sqrt{ \var { U^{i,l(i,v)}_{j,k(j,u)}  / {\cal F}^i_{j+1} }}               
      \left( \frac{ \gamma 2^{i+j-N} }{4} \right)^{-1/2} -1 \right)
      V^{i,l(i,v)}_{j,k(j,u)} \\
  + U^{i+1,l(i,v)/2}_{j+1,k(j,u)/2} \; 
  \delta^{i+1,l(i,v)/2}_{j,k(j,u)} \; \tilde{\delta}^{i,l(i,v)}_{j+1,k(j,u)/2}.
\end{eqnarray*}
The last term is equal to
$  \esp { < M|\tilde{e}_{l(i,v)}^i \otimes \tilde{e}_{k(j,u)}^j>
/ {\cal F}^i_{j+1} }$, it should not exist and its control is not straight.
Let us define the variables 
\( \X,\; i=B^*,\ldots,N-1,j+M(i)+1,\ldots,N-1 \)  by 
\begin{equation} \label{norm} \X=
 \left( \frac{ \gamma 2^{i+j-N} }{4} \right)^{-1/2} 
      V^{i,l(i,v)}_{j,k(j,u)}.
\end{equation}
The crucial point is that the variable $\X$
is ${\cal F}^i_{j}$ measurable, has ${\cal N}(0,1)$ distribution,
and is independent of ${\cal F}^i_{j+1}$. 
In particular, the variables 
\( \X,\; i=B^*,\ldots,N-1,j+M(i)+1,\ldots,N-1 \) 
are mutually independent.\\
By setting  
\begin{eqnarray*}
{\Delta_A}^{i,l(i,v)}_{j,k(j,u)}= 
c^i_v c^j_u \left( U^{i,l(i,v)}_{j,k(j,u)}-
  \esp { U^{i,l(i,v)}_{j,k(j,u)} / {\cal F}^i_{j+1} } 
  - \sqrt{ \var { U^{i,l(i,v)}_{j,k(j,u)}  / {\cal F}^i_{j+1} } }               
     \X \right), 
\\
   {\Delta_B}^{i,l(i,v)}_{j,k(j,u)}=c^i_v c^j_u
    \left(  \sqrt{ \var { U^{i,l(i,v)}_{j,k(j,u)}  / {\cal F}^i_{j+1} }  }             
      -\sqrt{ \left( \frac{ \gamma 2^{i+j-N} }{4} \right) } \right) 
     \X, 
\\
   {\Delta_C}^{i,l(i,v)}_{j,k(j,u)}= 
   \frac{c^i_v c^j_u}{4} U^{i+1,l(i,v)/2}_{j+1,k(j,u)/2} \; 
  \delta^{i+1,l(i,v)/2}_{j,k(j,u)} \; \tilde{\delta}^{i,l(i,v)}_{j+1,k(j,u)/2},
\end{eqnarray*} 
and, for $D \in \{ A,B,C \}$, 
\begin{eqnarray*} 
  T_2^D(u,v)= \sum_{i=B^*}^{N-1} \sum_{j=M(i)+1}^{N-1} 
   {\Delta_D}^{i,l(i,v)}_{j,k(j,u)} ,
   \\
   Q^D(u,v)= \proba{ \left\{
   |T_2^D(u,v)| \geq \frac{(x+\T\log (nab))^{3/2}}{24}
   \right\} \cap \Theta(u,v) },
\end{eqnarray*}
we obtain :
\begin{equation*} 
P_2(u,v) \leq Q^A(u,v)+Q^B(u,v)+Q^C(u,v).
\end{equation*}
\subsection{Control of  $Q^A(u,v)$. } \label{derA}
On $\Theta_0(u,v)$, we have
\[  \left| \delta^{i,l(i,v)}_{j,k(j,u)} \right| \leq \epsilon \mbox{ and }
\left| \tilde{\delta}^{i,l(i,v)}_{j,k(j,u)}  \right|
\leq \epsilon,
\]
and thus we can apply Lemma \ref{hyper} :
\begin{equation}
| \DA    | \leq c^i_v c^j_u \left( \alpha +\beta |\X |^2 \right) .
\label{cruc} 
\end{equation}
First we verify the conditions of Theorem \ref{penafin}. 
The sequence
\begin{eqnarray*}  {\Delta_A}^{N-1,l(N-1,v)}_{N-1,k(N-1,u)},
\ldots,
{\Delta_A}^{N-1,l(N-1,v)}_{M(N-1)+1,k(M(N-1)+1,u)},
{\Delta_A}^{N-2,l(N-2,v)}_{N-1,k(N-1,u)}, 
\ldots, \\ \\
{\Delta_A}^{N-2,l(N-2,v)}_{M(N-1)+1,k(M(N-1)+1,u)}, 
\ldots,  {\Delta_A}^{B^*,l(B^*,v)}_{A-1,k(A-1,u)} 
\end{eqnarray*}  
is adapted to the decreasing filtration
\[ \F^{N-1}_{N-1} \subset \ldots \subset
    \F^{N-1}_{M(N-1)+1} \subset \F^{N-2}_{N-1} \subset \ldots \subset
    \F^{N-2}_{M(N-2)+1} \subset \ldots \subset \F^{B^*}_{A-1}
\] 
because the variable $\DA$ is $\F^i_j$ measurable. Moreover,
\( \esp { \DA    / \F^i_{j+1} }=0 \). As in Theorem \ref{pena}, let
\[ \SI=\esp { \left( \DA    \right)^2 / \F^i_{j+1} }
    \mbox{ and } 
    {{\cal V}^{B^*}_{A-1}}^2 (u,v) = \sum_{i=B^*}^{N-1} \sum_{j=M(i)+1}^{N-1}  
    \SI .
\]
Using  (\ref{cruc}) and the properties of $\X$
(see (\ref{norm}) and its comment), one bounds $\SI$ :
\[ \SI \leq 2(c^i_vc^j_u)^2 (\alpha^2+3\beta^2), \]
and this allows us to bound ${{\cal V}^{B^*}_{A-1}}^2 (u,v)$. 
We obtain
\begin{eqnarray*} 
{{\cal V}^{B^*}_{A-1}}^2 (u,v) \leq \sum_{i=B^*}^{B-1} \sum_{j=M(i)+1}^{A-1} 
\frac{(\alpha^2+3\beta^2)}{8}+
\sum_{i=B^*}^{B-1} \sum_{j=A}^{N-1} 
\frac{(\alpha^2+3\beta^2)}{8} \left( \frac{2^A}{2^j} \right)^2 
\sum_{i=B}^{N-1} \sum_{j=A^*}^{A-1} 
\frac{(\alpha^2+3\beta^2)}{8} \left( \frac{2^B}{2^i} \right)^2 \\ +
\sum_{i=B}^{N-1} \sum_{j=A}^{N-1} 
\frac{(\alpha^2+3\beta^2)}{8} \left( \frac{2^A}{2^j} \right)^2
\left( \frac{2^B}{2^i} \right)^2
\end{eqnarray*}
with the convention $\sum_{s=N}^{N-1}=0$. It comes
\begin{eqnarray*} 
{{\cal V}^{B^*}_{A-1}}^2 (u,v) \leq \frac{(\alpha^2+3\beta^2)}{8}
\left( \frac{(B-B^*)(B-B^*+1)}{2}+\frac{4}{3}(B-B^*)
+\frac{4}{3}(A-A^*)+\frac{16}{9} \right) \\
\leq \theta (x+\T \log(nab) )^2,
\end{eqnarray*}
with $\theta = (\alpha^2+3\beta^2)/(8 (\T \ln (2))^2$. 
Moreover Inequality (\ref{cruc}) and Lemma \ref{coeff} give 
(\ref{cond2}) ; remark that it is only here that (\ref{cond3})
is not available. We obtain
\begin{eqnarray*}
    \esp { \left| \DA  \right|^{k}  /\F^i_{j+1} }
    \leq  \esp { 2^{k-1} (c^i_v c^j_u )^{k} 
    \left( \alpha^{k} +\beta^{k} \left| \X \right|^{2k} \right) /\F^i_{j+1} }
\\
     \leq \frac{1}{2} \left( \frac{\alpha}{2} \right)^{k}
    +\frac{1}{2} \left( \frac{\beta}{2} \right)^{k}
    \frac{(2k)!}{2^{k} k!} 
\\ 
    \leq \frac{1}{2} \left( \frac{\alpha}{2} \right)^{k}
    +\frac{1}{2} \left( \frac{\beta}{2} \right)^{k}
    2^k\frac{k!}{2} 
\\
    \leq   \frac{k!}{2} c^k
\end{eqnarray*}
with \[ c=\max ( \frac{\alpha}{2} , \beta) . \]
We can now apply Theorem \ref{penafin} :
\begin{eqnarray*}
Q^A(u,v)
&\leq & 2 \exp \left( \frac{-((x/2)+\T\log (nab))^{3}/{24}^2}
{2\left( \theta ((x/2)+\T \log(nab) )^2
+c((x/2)+\T\log (nab))^{3/2}/{24} \right) } \right) 
\\
&\leq & R_4 \exp(-\gamma_4 ((x/2)+\T\log (nab))) 
\end{eqnarray*}
where $R_4,\gamma_4$ are absolute positive constants
(we use $\T \geq 10$ and  $nab >496$ thus constants
do not depend on $\T$). In order to get
\[ Q^A(u,v) \leq R_4 \exp(-\gamma_4 (x/2) -2 \log(nab)) \]
we have to impose $\T \geq 2/\gamma_4$. 
\subsection{Control of  $Q^B(u,v)$. } \label{derB}
First we verify the conditions of Theorem \ref{penafin}. 
The sequence
\begin{eqnarray*}  {\Delta_B}^{N-1,l(N-1,v)}_{N-1,k(N-1,u)},
\ldots,
{\Delta_B}^{N-1,l(N-1,v)}_{M(N-1)+1,k(M(N-1)+1,u)},
{\Delta_B}^{N-2,l(N-2,v)}_{N-1,k(N-1,u)}, 
\ldots, \\ \\
{\Delta_B}^{N-2,l(N-2,v)}_{M(N-1)+1,k(M(N-1)+1,u)}, 
\ldots,  {\Delta_B}^{B^*,l(B^*,v)}_{A-1,k(A-1,u)} 
\end{eqnarray*}  
is adapted to the decreasing filtration
\[ \F^{N-1}_{N-1} \subset \ldots \subset
    \F^{N-1}_{M(N-1)+1} \subset \F^{N-2}_{N-1} \subset \ldots \subset
    \F^{N-2}_{M(N-2)+1} \subset \ldots \subset \F^{B^*}_{A-1}
\]
because the variable $\DB$ is $\F^i_j$ measurable. Moreover,
\( \esp { \DB   / \F^i_{j+1} }=0 \). 
As in Theorem \ref{pena}, let
\[ \SI=\esp { \left( \DB    \right)^2 / \F^i_{j+1} }
    \mbox{ et } 
    {{\cal V}^{B^*}_{A-1}}^2 (u,v) = \sum_{i=B^*}^{N-1} \sum_{j=M(i)+1}^{N-1}  
    \SI .
\]
The control of $Q^B(u,v)$ is based on the following lemma,
that we will prove in Section \ref{derB2}.
\begin{lem} \label{varq2b} On $\Theta(u,v)$, one gets : \\ a)
\[  (c^i_v c^j_u)^2
  \frac{   \left(  \var { U^{i,l(i,v)}_{j,k(j,u)}  / {\cal F}^i_{j+1} }               
   -\left(
   \displaystyle{ \frac{ \gamma 2^{i+j-N} }{4} }  \right) \right)^2 }
   {\gamma 2^{i+j-N+2} } \leq \theta_a ((x/2)+\T \log (nab))  \]
and \\ b)
\[ \sum_{i=B^*}^{N-1} \sum_{j=M(i)+1}^{N-1}
   (c^i_v c^j_u)^2
   \frac{   \left(  \var { U^{i,l(i,v)}_{j,k(j,u)}  / {\cal F}^i_{j+1} }               
    -\left(
    \displaystyle{ \frac{ \gamma 2^{i+j-N} }{4} }  \right) \right)^2 }
   {\gamma 2^{i+j-N+2} } \leq \theta_b ((x/2)+\T \log (nab))^{2} \]
where $\theta_a,\theta_b$ are absolute positive constants.
\end{lem}
\subsubsection{End of the control of $Q^B(u,v)$.}
Using the properties of $\X$ (see (\ref{norm}) and its comment), 
one gets on $\Theta(u,v)$
\begin{eqnarray*} \SI= (c^i_v c^j_u)^2
    \left(  \sqrt{ \var { U^{i,l(i,v)}_{j,k(j,u)}  / {\cal F}^i_{j+1} }  }             
     -\sqrt{ \left( \frac{ \gamma 2^{i+j-N} }{4} \right) } \right)^2
\\ \leq 4
  (c^i_v c^j_u)^2
   \frac{   \left(  \var { U^{i,l(i,v)}_{j,k(j,u)}  / {\cal F}^i_{j+1} }               
    -\left(
    \displaystyle{ \frac{ \gamma 2^{i+j-N} }{4} }  \right) \right)^2 }
   {\gamma 2^{i+j-N+2} }.
\end{eqnarray*}
Lemma \ref{varq2b} directly provides the bound of
${{\cal V}^{B^*}_{A-1}}^2 (u,v)$ on $\Theta(u,v)$ :
\[ {{\cal V}^{B^*}_{A-1}}^2 (u,v) \leq 4 \theta_b ((x/2)+\T \log (nab))^{2}
\]
and allows us to verify condition (\ref{cond3}) : on $\Theta(u,v)$,
using again the properties of $\X$ (see (\ref{norm}) and its comment),
one gets
\begin{eqnarray*}
\esp { \left( \DB \right)^{2k}  / \F^i_{j+1} }= 
 (c^i_v c^j_u)^{2k}
    \left(  \sqrt{ \var { U^{i,l(i,v)}_{j,k(j,u)}  / {\cal F}^i_{j+1} }  }             
     -\sqrt{ \left( \frac{ \gamma 2^{i+j-N} }{4} \right) } \right)^{2k}
     \frac{(2k)!}{k!2^k}
\\ \leq 
  (c^i_v c^j_u)^{2k} 4^k 
  \left(  \frac{   \left(  \var { U^{i,l(i,v)}_{j,k(j,u)}  / {\cal F}^i_{j+1} }               
    -\left(
    \displaystyle{ \frac{ \gamma 2^{i+j-N} }{4} }  \right) \right)^2 }
   {\gamma 2^{i+j-N+2} } \right)^k  \frac{(2k)!}{k!2^k}
\\ \leq  \frac{(2k)!}{k!2^k} c^{2k}
\end{eqnarray*}
with \[ c= \left( 4 \theta_a ((x/2)+\T \log (nab)) \right)^{1/2}. \]
We can now apply Theorem \ref{penafin}:
\begin{eqnarray*}
Q^B(u,v)
&\leq & 2 \exp \left( \frac{-((x/2)+\T\log (nab))^{3}/{24}^2}
{2\left( 4\theta_b ((x/2)+\T \log(nab) )^2
+2\sqrt{\theta_a}((x/2)+\T\log (nab))^{2}/{24} \right) } \right) 
\\
&\leq & R_5 \exp(-\gamma_5 ((x/2)+\T\log (nab))) 
\end{eqnarray*}
where $R_5,\gamma_5$ are absolute positive constants
(we use $\T \geq 10$ and $nab >496$, thus constants do not depend on
$\T$). In order to obtain
\[ Q^B(u,v) \leq R_5 \exp(-\gamma_5 (x/2) -2 \log(nab)) \]
we have to impose $\T \geq 2/\gamma_5$. 
\subsubsection{Proof of Lemma \ref{varq2b}.} \label{derB2}
We write 
\(   \displaystyle{ \var { U^{i,l(i,v)}_{j,k(j,u)}  / {\cal F}^i_{j+1} }               
      - \gamma \frac{ 2^{i+j-N+2} }{16}    } \) as a sum of three terms :
\begin{eqnarray*} 
\var { U^{i,l(i,v)}_{j,k(j,u)}  / {\cal F}^i_{j+1} }               
      &-& \gamma \frac{ 2^{i+j-N+2} }{16} 
= \frac{1}{16} \left( U^{i+1,l(i,v)/2}_{j+1,k(j,u)/2} 
- \gamma 2^{i+j-N+2} \right)
\\
&+& \frac{ U^{i+1,l(i,v)/2}_{j+1,k(j,u)/2} }{4}
\left( \frac{ U^{i+1,l(i,v)/2}_{j,k(j,u)} U^{i+1,l(i,v)/2}_{j,k(j,u)+1} }
{ \left( U^{i+1,l(i,v)/2}_{j+1,k(j,u)/2} \right)^2 }-\frac{1}{4} \right)
\\
&+& U^{i+1,l(i,v)/2}_{j+1,k(j,u)/2} 
\left( \frac{ U^{i+1,l(i,v)/2}_{j,k(j,u)} U^{i+1,l(i,v)/2}_{j,k(j,u)+1} }
{ \left( U^{i+1,l(i,v)/2}_{j+1,k(j,u)/2} \right)^2 } \right)
\left( \frac{ U^{i,l(i,v)}_{j+1,k(j,u)/2} U^{i,l(i,v)+1}_{j+1,k(j,u)/2} }
{ \left( U^{i+1,l(i,v)/2}_{j+1,k(j,u)/2} \right)^2 }-\frac{1}{4} \right).
\end{eqnarray*}
>From there, using $(a+b+c)^2 \leq 3(a^2+b^2+c^2)$ and, for the two last
terms, the relation 
\[ \frac{A}{T}(1-\frac{A}{T})-\frac{1}{4}=-\frac{(2A-T)^2}{4T^2}, \]
using moreover the properties of coefficients (see Lemma \ref{coeff})
and the following notation already used in the Section \ref{5} :
\[  \Delta^{i+1,l(i,v)/2}_{j,k(j,u)}= ( \alpha_j \beta_i )^{1/2}
                   \frac{ \left( 2U^{i+1,l(i,v)/2}_{j,k(j,u)}
                   -U^{i+1,l(i,v)/2}_{j+1,k(j,u)/2} \right)^2 }
                   {U^{i+1,l(i,v)/2}_{j+1,k(j,u)/2}} ,                  
\]
\[                  
      \tilde{\Delta}^{i,l(i,v)}_{j+1,k(j,u)/2}=(\alpha_j \beta_i  )^{1/2}
      \frac{ \left( 2 U^{i,l(i,v)}_{j+1,k(j,u)/2}-U^{i+1,l(i,v)/2}
      _{j+1,k(j,u)/2} \right)^2 }{U^{i+1,l(i,v)/2}_{j+1,k(j,u)/2}},
\]
we obtain :
\begin{eqnarray}
  (c^i_v c^j_u)
  \frac{   \left(  \var { U^{i,l(i,v)}_{j,k(j,u)}  / {\cal F}^i_{j+1} }               
   -\left(
   \displaystyle{ \frac{ \gamma 2^{i+j-N} }{4} }  \right) \right)^2 }
   {\gamma 2^{i+j-N+2} } \leq \frac{3}{2^8} (c^i_v c^j_u)
   \frac{ \left( U^{i+1,l(i,v)/2}_{j+1,k(j,u)/2} - \gamma 2^{i+j-N+2} \right)^2 }
   {\gamma 2^{i+j-N+2} } \nonumber
\\ \nonumber \\
+\frac{3}{2^8} \left( 
\frac{  
\left( \Delta^{i+1,l(i,v)/2}_{j,k(j,u)} \right)^2 + 
\left( \tilde{\Delta}^{i,l(i,v)}_{j+1,k(j,u)/2} \right)^2 }
{\gamma 2^{i+j-N+2} } \right)  . \label{rel1}
\end{eqnarray}
{\bf Bound of  } \( (c^i_v c^j_u)
   \displaystyle{ \frac{ \left( U^{i+1,l(i,v)/2}_{j+1,k(j,u)/2} - 
           \gamma 2^{i+j-N+2} \right)^2 }
   {\gamma 2^{i+j-N+2} } } \). \\ \\
With the convention $\sum_{s=N}^{N-1}=0$, the expansion of 
\( {U^{i+1,l(i,v)/2}_{j+1,k(j,u)/2}}-{\gamma 2^{i+j-N+2} } \)
on the basis ${\cal B}$ (defined by Section \ref{4}) is
\begin{eqnarray*} 
{U^{i+1,l(i,v)/2}_{j+1,k(j,u)/2}}-{\gamma 2^{i+j-N+2} }
=<M|e^{i+1}_{l(i,v)/2}\otimes e^{j+1}_{k(j,u)/2}>-n\frac{2^{i+1}}{2^N}
\frac{2^{j+1}}{2^N} <M|e^{N}_{0}\otimes e^{N}_{0}> \\
=\sum_{s=i+1}^{N-1} \pm \frac{2^{i+1}}{2^{s+1}}
<M|\tilde{e}^{s}_{l(s,v)/2}\otimes e^{j+1}_{k(j,u)/2}>
+\frac{2^{i+1}}{2^N}
\left( <M|e^{i+1}_{l(i,v)/2}\otimes e^{j+1}_{k(j,u)/2}> - \frac{2^{j+1}}{2^N}
<M|e^{N}_{0}\otimes e^{N}_{0}> \right) \\
=\sum_{s=i+1}^{N-1} \pm \frac{2^{i+1}}{2^{s+1}}
<M|\tilde{e}^{s}_{l(s,v)/2} \otimes e^{j+1}_{k(j,u)/2}>
+\frac{2^{i+1}}{2^N} \sum_{r=j+1}^{N-1} \pm \frac{2^{j+1}}{2^{r+1}}
<M|{{e}}^{N}_{0} \otimes {\tilde{e}}^{r}_{k(r,u)}> .
\end{eqnarray*}
Let us recall that on $\Theta_0(u,v)$ we have
\[ \left\{ \begin{array}{l}
{U^{s+1,l(s,v)/2}_{j+1,k(j,u)/2}} \leq \gamma (1+\epsilon) 2^{s+j+2-N} \\ \\
{U^{N,0}_{r+1,k(r,v)}} \leq \gamma (1+\epsilon) 2^{r+1-N} .
\end{array} \right.
\]
This yields
\begin{eqnarray*}
(c^i_v c^j_u)
   \displaystyle{ \frac{ \left( U^{i+1,l(i,v)/2}_{j+1,k(j,u)/2} - 
           \gamma 2^{i+j-N+2} \right)^2 }
   {\gamma 2^{i+j-N+2} } } \\
   \leq 2(1+\epsilon) (c^i_v c^j_u)^{1/2}  \left[
   \left( \sum_{s=i+1}^{N-1} \left( \frac{2^{i+1}}{2^{s+1}} \right)^{1/4}
   (c^i_v c^j_u)^{1/4} \left( \frac{2^{i+1}}{2^{s+1}} \right)^{3/4}
   \frac{2^{(s+j+2-N)/2}}{2^{(i+j+2-N)/2}}
   \frac{<M|\tilde{e}^{s}_{l(s,v)/2}\otimes e^{j+1}_{k(j,u)/2}>}{
   \sqrt{ {U^{s+1,l(s,v)/2}_{j+1,k(j,u)/2}} }} \right)^2 \right.
   \\
   \left. +
   \left( \sum_{r=j+1}^{N-1} \left( \frac{2^{j+1}}{2^{r+1}} \right)^{1/4}
   \left( \frac{2^{i+1}}{2^{N}} \right)^{1/4}
   (c^i_v c^j_u)^{1/4} \left( \frac{2^{j+1}}{2^{r+1}} \right)^{3/4}
   \left( \frac{2^{i+1}}{2^{N}} \right)^{3/4}
   \frac{2^{(r+1)/2}}{2^{(i+j+2-N)/2}}
   \frac{<M|{e}^{N}_{0}\otimes e^{r}_{k(r,u)/2}>}{
   \sqrt{ {U^{N,0}_{r+1,k(r,u)/2}} }} \right)^2 \right].
\end{eqnarray*}
The relations
\[ <M|\tilde{e}^{s}_{l(s,v)/2}\otimes e^{j+1}_{k(j,u)/2}>=
    \frac{ {U^{s,l(s,v)}_{j+1,k(j,u)/2}}-
    {U^{s,l(s,v)+1}_{j+1,k(j,u)/2}} }{ {U^{s+1,l(s,v)/2}_{j+1,k(j,u)/2}} }
\]
and
\[ <M|{e}^{N}_{0}\otimes \tilde{e}^{r}_{k(r,v)}>=
    \frac{ {U^{N,0}_{r,k(r,v)}}-{U^{N,0}_{r,k(r,v)+1}} }
    { {U^{N,0}_{r+1,k(r,v)/2}} } ,
\]
as well as the relations (see Lemma \ref{coeff})
\[ \frac{2^{i+1}}{2^{s+1}} c^i_v \leq \beta_s,\;
    \frac{2^{j+1}}{2^{r+1}} c^j_u \leq \alpha_r ,
\]
give us the following inequality
\begin{eqnarray} (c^i_v c^j_u)
   \displaystyle{ \frac{ \left( U^{i+1,l(i,v)/2}_{j+1,k(j,u)/2} - 
           \gamma 2^{i+j-N+2} \right)^2 }
   {\gamma 2^{i+j-N+2} } } 
   \leq 2(1+\epsilon) (c^i_v c^j_u)^{1/2}  \left[
   \left( \sum_{s=i+1}^{N-1} \left( \frac{2^{i}}{2^{s}} \right)^{1/4}
   \left( \tilde{\Delta}^{s,l(s,v)}_{j+1,k(j,u)/2} \right)^{1/2}
   \right)^{2} \right. \nonumber \\
   + \left. \left( \sum_{r=j+1}^{N-1} \left( \frac{2^{j}}{2^{r}} \right)^{1/4}
   \left( \frac{2^{i+1}}{2^{N}} \right)^{1/4}
   \left( {\Delta}^{N,0}_{r,k(r,u)} \right)^{1/2}
   \right)^{2} \right] . \label{rel2}
\end{eqnarray}
{\bf Proof of Inequality a) of Lemma \ref{varq2b}.} \\ \\
On  $\Theta_0(u,v)$ one gets
\begin{eqnarray} \frac{ \Delta^{i+1,l(i,v)/2}_{j,k(j,u)} }{\gamma 2^{i+j+2-N}}
    = (\beta_i \alpha_i)^{1/2} 
    \frac{ U^{i+1,l(i,v)/2}_{j+1,k(j,u)/2} }{\gamma 2^{i+j+2-N}}
    \left( \delta^{i+1,l(i,v)/2}_{j,k(j,u)} \right)^2
    \leq \epsilon^2 (1+\epsilon)/2 \nonumber \\ \label{rel3} \\
    \mbox{ and }  \;\;
    \frac{ \tilde{\Delta}^{i,l(i,v)}_{j+1,k(j,u)/2} }{\gamma 2^{i+j+2-N}}
    = (\beta_i \alpha_i)^{1/2} 
    \frac{ U^{i+1,l(i,v)/2}_{j+1,k(j,u)/2} }{\gamma 2^{i+j+2-N}}
    \left( \tilde{\delta}^{i+1,l(i,v)/2}_{j,k(j,u)} \right)^2
    \leq \epsilon^2 (1+\epsilon)/2. \nonumber
\end{eqnarray}
Using Cauchy-Schwarz Inequality, as well as Inequalities
(\ref{rel1}), (\ref{rel2}), (\ref{rel3}), we obtain on $\Theta(u,v)$
\begin{eqnarray*}
(c^i_v c^j_u)^2
   \displaystyle{ \frac{ \left( U^{i+1,l(i,v)/2}_{j+1,k(j,u)/2} - 
           \gamma 2^{i+j-N+2} \right)^2 }
   {\gamma 2^{i+j-N+2} } } \leq \frac{12(1+\epsilon)}{2^8}
   (c^i_v c^j_u)^{3/2} \frac{1}{\sqrt{2}-1} ((x/2)+\T \log (nab)) \\
   +\frac{3\epsilon^2 (1+\epsilon) }{2^8} (c^i_v c^j_u)
   \left( \Delta^{i+1,l(i,v)/2}_{j,k(j,u)} +
           \tilde{\Delta}^{i,l(i,v)}_{j+1,k(j,u)/2} \right) \\
    \leq \theta_a ((x/2)+\T \log (nab))
\end{eqnarray*}
with
\[ \theta_a=\frac{3 (1+\epsilon) }{2^9} \left( \frac{1}{\sqrt{2}-1}
    +\epsilon^2/2 \right).
\]
{\bf Proof of Inequality b) of Lemma \ref{varq2b}.} \\ \\
Using Cauchy-Schwarz Inequality, as well as Inequalities
(\ref{rel1}), (\ref{rel2}), (\ref{rel3}), we obtain on $\Theta(u,v)$
\begin{eqnarray*} 
\sum_{i=B^*}^{N-1} \sum_{j=M(i)+1}^{N-1}
   (c^i_v c^j_u)^2
   \frac{   \left(  \var { U^{i,l(i,v)}_{j,k(j,u)}  / {\cal F}^i_{j+1} }               
    -\left(
    \displaystyle{ \frac{ \gamma 2^{i+j-N} }{4} }  \right) \right)^2 }
   {\gamma 2^{i+j-N+2} } 
   \\ \leq 
    \frac{6(1+\epsilon)}{2^8} \left( \frac{1}{2^{1/4}-1} \right)
    \sum_{i=B^*}^{N-1} \sum_{j=M(i)+1}^{N-1}
    (c^i_v c^j_u)^{3/2} \sum_{s=i+1}^{N-1} 
    \left( \frac{2^{i}}{2^{s}} \right)^{1/4}
    \tilde{\Delta}^{s,l(s,v)}_{j+1,k(j,u)/2} 
    \\
   + \frac{6(1+\epsilon)}{2^8} \left( \frac{1}{2^{1/4}-1} \right)
    \sum_{i=B^*}^{N-1} \sum_{j=M(i)+1}^{N-1} (c^i_v c^j_u)^{3/2}
    \sum_{r=j+1}^{N-1} \left( \frac{2^{j}}{2^{r}} \right)^{1/4}
    {\Delta}^{N,0}_{r,k(r,u)} \\
   + \frac{3 \epsilon^2 (1+\epsilon)}{2^8} 
   \sum_{i=B^*}^{N-1} \sum_{j=M(i)+1}^{N-1} (c^i_v c^j_u)
   \left( \Delta^{i+1,l(i,v)/2}_{j,k(j,u)} +
           \tilde{\Delta}^{i,l(i,v)}_{j+1,k(j,u)/2} \right).
\end{eqnarray*}
Exchanging the sums and considering the definitions of $M(i)$, 
${\cal M}(j)$, one gets on $\Theta(u,v)$ :
\begin{eqnarray*} 
\sum_{i=B^*}^{N-1} \sum_{j=M(i)+1}^{N-1}
   (c^i_v c^j_u)^2
   \frac{   \left(  \var { U^{i,l(i,v)}_{j,k(j,u)}  / {\cal F}^i_{j+1} }               
    -\left(
    \displaystyle{ \frac{ \gamma 2^{i+j-N} }{4} }  \right) \right)^2 }
   {\gamma 2^{i+j-N+2} } 
   \\ \leq 
    \frac{6(1+\epsilon)}{2^8} \left( \frac{1}{2^{1/4}-1} \right)^2
    \left( \frac{1}{2} \right)^2
    \left( \sum_{j=A^*}^{N-1} c^j_u + \sum_{i=B^*}^{N-1} c^i_u \right)
    ((x/2)+\T \log (nab))   
    \\ + \frac{3 \epsilon^2 (1+\epsilon)}{2^9} 
    \left( \sum_{i=B^*}^{N-1} c^i_u + \sum_{j=A^*}^{N-1} c^j_u \right)
    ((x/2)+\T \log (nab)).
\end{eqnarray*}
Then using the properties of coefficients 
(see Lemma \ref{coeff}) we obtain on $\Theta(u,v)$
\[ \sum_{i=B^*}^{N-1} \sum_{j=M(i)+1}^{N-1}
   (c^i_v c^j_u)^2
   \frac{   \left(  \var { U^{i,l(i,v)}_{j,k(j,u)}  / {\cal F}^i_{j+1} }               
    -\left(
    \displaystyle{ \frac{ \gamma 2^{i+j-N} }{4} }  \right) \right)^2 }
   {\gamma 2^{i+j-N+2} } \leq  \theta_b ((x/2)+\T \log (nab))^2
\]
with \[ \theta_b= \frac{3 (1+\epsilon)}{2^9} 
         \left( \left( \frac{1}{2^{1/4}-1} \right)^2 + \epsilon^2 \right)
         \left( 1+\frac{1}{10\log(2)} \right).
      \]
\subsection{Control of $Q^C(u,v)$.}
As mentioned previously, this term comes, to our mind, from a flaw of the
construction. Its control is closer to the control of
$P_1^E(u,v)$ (Section \ref{62}) than to the controls of
$Q^A(u,v)$, $Q^B(u,v)$ (Sections \ref{derA} and \ref{derB}). The variable
$\DC$ is ${\cal F}^i_{j+1}$ measurable, but
\( \esp { \DC/{\cal F}^i_{j+2} } \neq 0 \). Therefore, in order to apply
Theorem \ref{penafin}, we have to consider the variables $\DDC$
defined by
\[ \DDC=\sum_{j=M(i)+1}^{N-1} \frac{c^i_v c^j_u}{4} 
   U^{i+1,l(i,v)/2}_{j+1,k(j,u)/2} \; 
  \delta^{i+1,l(i,v)/2}_{j,k(j,u)} \; 
  \tilde{\delta}^{i,l(i,v)}_{j+1,k(j,u)/2}.
\]
Thus the term to be controlled is equal to 
\[ T_2^C(u,v)=\sum_{i=B^*}^{N-1} \DDC. \]
We verify the conditions of Theorem \ref{penafin}. 
The sequence
\(  \DDC \; ,i=N-1,\ldots,B^* \)
is adapted to the decreasing filtration 
\( \F^{N-1}_{0} \subset
    \F^{N-2}_{0} \subset  \ldots \subset \F^{B^*}_{0}
\)
because the variable $\DDC$ is $\F^i_0$ measurable. 
In order to verify the other conditions of Theorem \ref{penafin}
easily, we start by giving a new expression for $\DDC$.
\subsubsection{A new expression for $\DDC$.}
With the notations of Section \ref{2}, one gets
\[ \DDC=\sum_{j=M(i)+1}^{N-1} \frac{c^i_v c^j_u}{4} 
  \delta^{i+1,l(i,v)/2}_{j,k(j,u)} \; 
  <M|\tilde{e}^{i}_{l(i,v)}\otimes e^{j+1}_{k(j,u)/2}>.
\]
Let us recall that $k(j,u)$ is defined by the even integer such that
\[  u2^{A^*} \in ]k(j,u)2^{j}; ( k(j,u) +2) 2^{j} . \]
This yields
\[ k(j,u)=\left\{  \begin{array}{ll} 
                        \displaystyle{\frac{k(j-1,u)}{2}} & \mbox{ if }
                        \displaystyle{\frac{k(j-1,u)}{2}}\mbox{ even, } \\
                        \displaystyle{\frac{k(j-1,u)}{2}} - 1 & \mbox{ if }
                        \displaystyle{\frac{k(j-1,u)}{2}} \mbox{ odd. }
                   \end{array} \right.
\]
We define $k^*(j-1,u)$ by :
\[ k^*(j-1,u)=\left\{  \begin{array}{ll} k(j-1,u)+{2} & \mbox{ if }
                       \displaystyle{\frac{k(j-1,u)}{2}} \mbox{ even, } \\
                        k(j-1,u)-{2}  & \mbox{ if  }
                        \displaystyle{\frac{k(j-1,u)}{2}} \mbox{ odd. }
                   \end{array} \right.
\]
In this way we have
\[ e^{j+1}_{\frac{k(j,u)}{2}}=e^{j}_{\frac{k(j-1,u)}{2}}+e^{j}_{\frac{k^*(j-1,u)}{2}}. \]
In other words, if one interprets the vector
$e^{j+1}_{\frac{k(j,u)}{2}}$ as a representation of the length $2^{j+1}$
interval containing $u2^{A^*}$,
then the vector $e^{j}_{\frac{k(j-1,u)}{2}}$ is a 
representation of the length $2^{j}$ interval containing $u2^{A^*}$
and is one half (left or right) of the length $2^{j+1}$
interval containing $u2^{A^*}$ ;
the vector $e^{j}_{\frac{k^*(j-1,u)}{2}}$ represents
the other half of the length $2^{j+1}$
interval containing $u2^{A^*}$,
that is to say the one not containing $u2^{A^*}$.
Using this relation, one gets
\begin{eqnarray*}
<M|\tilde{e}^{i}_{l(i,v)}\otimes e^{j+1}_{\frac{k(j,u)}{2}}>=
<M|\tilde{e}^{i}_{l(i,v)}\otimes \left( e^{j+1}_{\frac{k(j,u)}{2}}-e^{j}_{\frac{k(j-1,u)}{2}}
+e^{j}_{\frac{k(j-1,u)}{2}} \right)> 
\\ = <M|\tilde{e}^{i}_{l(i,v)}\otimes e^{j}_{\frac{k^*(j-1,u)}{2}}>
+<M|\tilde{e}^{i}_{l(i,v)}\otimes e^{j}_{\frac{k(j-1,u)}{2}}> 
\\ = \left( \sum_{r=M(i)+1}^j <M|\tilde{e}^{i}_{l(i,v)}\otimes e^{r}_{\frac{k^*(r-1,u)}{2}}>
\right) + <M|\tilde{e}^{i}_{l(i,v)}\otimes e^{M(i)+1}_{\frac{k(M(i),u)}{2}}> .
\end{eqnarray*}
>From there we have a new expression for $\DDC$ :
\begin{eqnarray*} \DDC= \left( \sum_{r=M(i)+1}^{N-1} \alpha^{i+1,l(i,v)/2}_r
   <M|\tilde{e}^{i}_{l(i,v)}\otimes e^{r}_{\frac{k^*(r-1,u)}{2}}> \right)
   \\ +\alpha_{M(i)+1} 
   <M|\tilde{e}^{i}_{l(i,v)}\otimes e^{M(i)+1}_{\frac{k(M(i),u)}{2}}>
\end{eqnarray*}
where $\alpha^{i+1,l(i,v)/2}_r$ is the random variable defined by
\[ \alpha^{i+1,l(i,v)/2}_r = \sum_{j=r}^{N-1} \frac{c^i_v c^j_u}{4} 
  \delta^{i+1,l(i,v)/2}_{j,k(j,u)}  .
\]
Remark that $\alpha^{i+1,l(i,v)/2}_r$ is ${\cal F}^{i+1}_0$ measurable.
The interest of the variables 
\[ <M|\tilde{e}^{i}_{l(i,v)}\otimes e^{M(i)+1}_{\frac{k(M(i),u)}{2}}>,
<M|\tilde{e}^{i}_{l(i,v)}\otimes e^{r}_{\frac{k^*(r-1,u)}{2}}>,\; 
r=M(i)+1,\ldots,N-1, \] 
is to be {\em nearly} independent, while the
variables \[ <M|\tilde{e}^{i}_{l(i,v)}\otimes e^{j+1}_{\frac{k(j,u)}{2}}>, \;
j=M(i)+1,\ldots,N-1, \] 
are closely correlated. More precisely, 
\begin{eqnarray} {\cal L} \left( <M|{e}^{i}_{l(i,v)}\otimes e^{M(i)+1}_{\frac{k(M(i),u)}{2}}>,
<M|{e}^{i}_{l(i,v)}\otimes e^{r}_{\frac{k^*(r-1,u)}{2}}>,\; 
r=M(i)+1,\ldots,N-1 
\right/ \left. {\cal F}^{i+1}_0 \right) =  \nonumber \\ 
{\cal B} 
\left( <M|{e}^{i+1}_{\frac{l(i,v)}{2}} \otimes e^{M(i)+1}_{\frac{k(M(i),u)}{2}} > ,1/2 \right)
\otimes \bigotimes_{r=M(i)+1}^{N-1}
{\cal B} 
\left( <M|{e}^{i+1}_{\frac{l(i,v)}{2}}\otimes e^{r}_{\frac{k^*(r-1,u)}{2}}> ,1/2 \right).
\label{prodloi}
\end{eqnarray} 
As in Section \ref{62} (Equality (\ref{bincent})) we have
\begin{eqnarray} \label{bnct}
<M|\tilde{e}^{i}_{l(i,v)}\otimes e^{r}_{\frac{k^*(r-1,u)}{2}}>
=2<M|{e}^{i}_{l(i,v)}\otimes e^{r}_{\frac{k^*(r-1,u)}{2}}>
-<M|{e}^{i+1}_{\frac{l(i,v)}{2}}\otimes e^{r}_{\frac{k^*(r-1,u)}{2}}>, \\
<M|\tilde{e}^{i}_{l(i,v)}\otimes e^{M(i)+1}_{\frac{k(M(i),u)}{2}}>=
2<M|{e}^{i}_{l(i,v)}\otimes e^{M(i)+1}_{\frac{k(M(i),u)}{2}}>
-<M|{e}^{i+1}_{\frac{l(i,v)}{2}}\otimes e^{M(i)+1}_{\frac{k(M(i),u)}{2}}>.
\nonumber
\end{eqnarray}
Using (\ref{prodloi}) and (\ref{bnct}) one gets
\begin{eqnarray*} {\cal L} \left( 
<M|\tilde{e}^{i}_{l(i,v)}\otimes e^{M(i)+1}_{\frac{k(M(i),u)}{2}}>,
<M|\tilde{e}^{i}_{l(i,v)}\otimes e^{r}_{\frac{k^*(r-1,u)}{2}}>,\; 
r=M(i)+1,\ldots,N-1 
\right/ \left. {\cal F}^{i+1}_0 \right) =  \\ 
{\cal L} \left( \sum_{u=1}^{M_1} X_u, \sum_{u=M_1+1}^{M_1+M_2} X_u,\ldots,
\sum_{u=M_1+\cdots+M_{N-M(i)-1}+1}^{M_1+\cdots+M_{N-M(i)}} X_u \right)
\end{eqnarray*}
where $M_s$, $s=1,\ldots,N-M(i)$, are the random variables defined by
\begin{eqnarray*}
M_1=<M|{e}^{i+1}_{\frac{l(i,v)}{2}}\otimes e^{M(i)+1}_{\frac{k(M(i),u)}{2}}>, \\
M_s=<M|{e}^{i+1}_{\frac{l(i,v)}{2}}\otimes e^{s+M(i)-1}_{\frac{k^*(s+M(i)-2,u)}{2}}>
\mbox{ pour } s=2,\ldots,N-M(i),
\end{eqnarray*}
and where $X_1,\ldots,X_{M_1+\cdots+M_{N-M(i)}}$ are i.i.d. random variables 
\( \proba { X_1=+1 } =\proba {X_1 = -1 }
=1/2 \). By setting 
\[
\beta_u=\left\{ \begin{array}{ll}
\alpha^{i+1,l(i,v)/2}_{M(i)+1} & \mbox{ for } u=1,\ldots,M_1, \\ & \\
\alpha^{i+1,l(i,v)/2}_{s+M(i)-1} & \mbox{ for } 
u=M_1+\cdots+M_{s-1}+1,\ldots,M_1+\cdots+M_{s}, s = 2, \ldots,N-M(i),
\end{array} \right.
\]
one obtains the very simple expression
\[ {\cal L} \left( \DDC \right/ \left. {\cal F}^{i+1}_0 \right) =
    \sum_{u=1}^{M_1+\cdots+M_{N-M(i)}} \beta_u X_u . 
\]
\subsubsection{End of the control of $Q_C(u,v)$.}
Clearly $\esp { \DDC / {\cal F}^{i+1}_0 } =0$.
As in Theorem \ref{pena}, let
\[ \left( \sigma^{i}(u,v)  \right)^2 = 
\esp { \left(  \DDC
\right)^2 / {\cal F}^{i+1}_0 } \mbox{ and }  
{\cal V}^2_{B^*}(u,v)=\sum_{i=B^*}^{N-1} \left( \sigma^{i}(u,v) \right)^2 . \]
We have
\[
    \left( \sigma^{i}(u,v)  \right)^2 = 
    \sum_{u=1}^{M_1+\cdots+M_{N-M(i)}} \beta_u^2  
    =M_1 \left( \alpha^{i+1,l(i,v)/2}_{M(i)+1} \right)^2
    +\sum_{s=2}^{N-M(i)} M_s \left( \alpha^{i+1,l(i,v)/2}_{s+M(i)-1} \right)^2.
\]
Moreover on $\Theta_0(u,v)$ we have
\begin{eqnarray*}
M_1 \leq \gamma (1+\epsilon) 2^{i+1+M(i)+1-N}, \\ \\
M_s \leq \gamma (1+\epsilon) 2^{i+1+s+M(i)-1-N} 
\mbox{ for } s=2,\ldots,N-M(i),
\end{eqnarray*}
this yields
\begin{eqnarray*} 
M_1 \left( \alpha^{i+1,l(i,v)/2}_{M(i)+1} \right)^2
&\leq & \gamma (1+\epsilon)
\left(   \sum_{j=M(i)+1}^{N-1} \frac{(c^i_v c^j_u)}{4} 
    \left| \delta^{i+1,l(i,v)/2}_{j,k(j,u)} \right| 
    \left( 2^{i+2+M(i)-N} \right)^{1/2} 
\right)^2
\\ \mbox{and for } s =2, \ldots, N-M(i), & \\
M_s \left( \alpha^{i+1,l(i,v)/2}_{s+M(i)-1} \right)^2
&\leq & \gamma (1+\epsilon)
\left(   \sum_{j=s+M(i)-1}^{N-1} \frac{(c^i_v c^j_u)}{4} 
    \left| \delta^{i+1,l(i,v)/2}_{j,k(j,u)} \right| 
    \left( 2^{i+s+M(i)-N} \right)^{1/2}   \right)^2 .
\end{eqnarray*}
Using the notation of Section \ref{5} :
\[  \Delta^{i+1,l(i,v)/2}_{j,k(j,u)}= (\alpha_j \beta_i )^{1/2}
                   \frac{ \left( U^{i+1,l(i,v)/2}_{j,k(j,u)}
                   -U^{i+1,l(i,v)/2}_{j,k(j,u)+1} \right)^2 }
                   {U^{i+1,l(i,v)/2}_{j+1,k(j,u)/2}}                   
\]
one gets on $\Theta_0(u,v)$ 
\begin{eqnarray*} 
M_1 \left( \alpha^{i+1,l(i,v)/2}_{M(i)+1} \right)^2
& \leq & \left( \sum_{j=M(i)+1}^{N-1} (c^i_v c^j_u)^{3/4} 
\left( \frac{2^{M(i)}}{{2}^{j}} \right)^{1/2}
\left( \Delta^{i+1,l(i,v)/2}_{j,k(j,u)} \right)^{1/2} \right)^2
\\
\mbox{and for } s =2, \ldots, N-M(i), & \\
M_s \left( \alpha^{i+1,l(i,v)/2}_{s+M(i)-1} \right)^2
& \leq & \frac{ (1+\epsilon)}{(1-\epsilon)}
\left( \sum_{j=s+M(i)-1}^{N-1} (c^i_v c^j_u)^{3/4} 
\left( \frac{2^{s+M(i)}}{{2}^{j+2}} \right)^{1/2}
\left( \Delta^{i+1,l(i,v)/2}_{j,k(j,u)} \right)^{1/2} \right)^2.
\end{eqnarray*}
With Cauchy-Schwarz Inequality we have :
\begin{eqnarray*} 
M_1 \left( \alpha^{i+1,l(i,v)/2}_{M(i)+1} \right)^2
&\leq & \frac{ (1+\epsilon)}{16(1-\epsilon)} 
\left( \frac{c^i_v}{2} \right)^{3/2}
\left( \frac{1}{\sqrt{2}-1} \right)  \sum_{j=M(i)+1}^{N-1} 
\left( \frac{2^{M(i)}}{{2}^{j}} \right)^{1/2}
\Delta^{i+1,l(i,v)/2}_{j,k(j,u)} 
\\
\mbox{and for } s =2, \ldots, N-M(i), & \\
M_s \left( \alpha^{i+1,l(i,v)/2}_{s+M(i)-1} \right)^2
& \leq & \frac{ (1+\epsilon)}{16(1-\epsilon)} 
\left( \frac{c^i_v}{2} \right)^{3/2}
\left( \frac{1}{\sqrt{2}-1} \right)  \sum_{j=s+M(i)-1}^{N-1} 
\left( \frac{2^{s+M(i)}}{{2}^{j+2}} \right)^{1/2}
\Delta^{i+1,l(i,v)/2}_{j,k(j,u)} .
\end{eqnarray*}
Now we exchange the sums, on $\Theta_1(u,v)$ this leads to :
\[ \left( \sigma^{i}(u,v)  \right)^2 \leq 
\frac{ (1+\epsilon)}{16(1-\epsilon)} \left( \frac{c^i_v}{2} \right)^{3/2}
\left( \frac{(2\sqrt{2})-1}{\left( \sqrt{2}-1 \right)^2} \right) 
\left( (x/2) +\T \log(nab) \right). 
\]
Then using the properties of coefficients
(see Lemma \ref{coeff})
and the convention $\sum_{i=N-1}^N=0$ we obtain :
\begin{eqnarray*} {\cal V}^2_{B^*}(u,v) \leq 
   \frac{ (1+\epsilon)}{16(1-\epsilon)} 
   \left( \frac{(2\sqrt{2})-1}{\left( \sqrt{2}-1 \right)^2} \right) 
   \left( (x/2) +\T \log(nab) \right) 
   \left( \sum_{i=B^*}^{B-1} \left( \frac{1}{2} \right)^{3/2}
   + \sum_{i=B}^{N-1} \left( \frac{2^B}{2^{i+1}} \right)^{3/2} \right)
   \\ \leq 0.03 \left( (x/2) +\T \log(nab) \right)^2.
\end{eqnarray*}
On the other hand, with the same argument as in the proof of Lemma
\ref{bineg}, one gets
\begin{eqnarray*} 
\esp { \left( \DDC \right)^{2k} \left/  {\cal F}^{i+1}_0 \right.  }
&\leq & \frac{ (2k)! }{k! 2^k} \sum_{u_1} \cdots \sum_{u_k}
\beta^2_{u_1} \ldots \beta^2_{u_k} \\ &= & \frac{ (2k)! }{k! 2^k} 
\left( \sum_{u} \beta^2_{u}\right)^k \\ &= & \frac{ (2k)! }{k! 2^k} 
\left( \sigma^{i}(u,v)  \right)^{2k} \\ &\leq & \frac{ (2k)! }{k! 2^k} c^{2k}
\end{eqnarray*}
with
\[ c = \left( \theta_c  (x/2) +\T \log(nab)  \right)^{1/2}
:=\left( \frac{ (1+\epsilon)}{(1-\epsilon)} 
\left( \frac{1}{2} \right)^{7}
\left( \frac{(2\sqrt{2})-1}{\left( \sqrt{2}-1 \right)^2} \right) 
\left( (x/2) +\T \log(nab) \right) \right)^{1/2}.
\]
We can now apply Theorem \ref{penafin}:
\begin{eqnarray*}
Q^C(u,v)
&\leq & 2 \exp \left( \frac{-((x/2)+\T\log (nab))^{3}/{24}^2}
{2\left( 0.03((x/2)+\T \log(nab) )^2
+\theta_c((x/2)+\T\log (nab))^{2}/{24} \right) } \right) 
\\
&\leq & 2 \exp \left(-\frac{1}{35}( \frac{x}{2}+\T\log (nab) ) \right) .
\end{eqnarray*}
In order to obtain
\[ Q^C(u,v) \leq 2 \exp(-(x/70) -2 \log(nab)) \]
we have to impose $\T \geq 70$. \\ \\


\end{document}